\documentclass[review]{elsarticle}
\usepackage[english]{babel}
\usepackage{a4wide}
\usepackage[ansinew]{inputenc}
\usepackage{amsfonts}
\usepackage{amsmath}
\usepackage{amssymb}
\usepackage{amsthm}
\usepackage{amstext}
\usepackage{newlfont}
\usepackage{graphicx}
\usepackage{color}
\usepackage{epstopdf}
\usepackage{here}
\usepackage{algorithm}
\usepackage{algpseudocode}

\newtheorem{theorem}{Theorem}[section]

\begin{document}

\begin{frontmatter}

\title{$M/G/c/c$ state dependent queuing model for a road traffic system of two sections in tandem}

\author[mymainaddress]{Nacira Guerouahane \corref{mycorrespondingauthor}}
\cortext[mycorrespondingauthor]{Corresponding author. Tel.: +213-34-21-08-00.}
\ead{naciraro@hotmail.fr}

\author[mymainaddress]{Djamil Aissani}
\ead{lamos\_bejaia@hotmail.com }

\author[mysecondaryaddress]{Nadir Farhi}
\ead{nadir.farhi@ifsttar.fr }

\author[mymainaddress]{Louiza Bouallouche-Medjkoune}
\ead{louiza\_medjkoune@yahoo.fr}

\address[mymainaddress]{LaMOS Research Unit, Faculty of Exact Sciences, University of Bejaia, 06000 Bejaia,  Algeria.}

\address[mysecondaryaddress]{Universit\'e Paris-Est, COSYS, GRETTIA, IFSTTAR, F-77447 Marne-la-Vallée, France.}

\begin{abstract}
We propose in this article a $M/G/c/c$ state dependent queuing model for road traffic flow.
The model is based on finite capacity queuing theory which captures the
stationary density-flow relationships. It is also
inspired from the deterministic Godunov scheme for the road traffic simulation. We first present a reformulation of
the existing linear case of $M/G/c/c$ state dependent model, in order to use flow rather than speed variables.
We then extend this model in order to consider upstream traffic demand and downstream traffic supply.
After that, we propose the model for two road sections in
tandem where both sections influence each other. In order to deal with this mutual dependence, we solve an implicit
system given by an algebraic equation.
Finally, we derive some performance measures (throughput and expected travel time). A comparison with results predicted
by the $M/G/c/c$ state dependent queuing networks shows that the model we propose here captures really the dynamics of
the road traffic.
\end{abstract}

\begin{keyword}
Traffic flow modeling \sep  finite queuing systems \sep $M/G/c/c$ state dependent.
\end{keyword}

\end{frontmatter}

\section{Introduction}

Congestion is a phenomenon that arises in a variety of contexts. The most familiar representation
is urban traffic congestion. Congested networks involve complex traffic interactions. Providing an
analytical description of these intricate interactions is challenging. The dynamics of traffic flows
is submitted to stochastic nature of traffic demand and supply functions of vehicles passing from one
section to another.

One of the well known macroscopic traffic model on a road section is the LWR
(Lighthill-Whitham-Richards) first-order continuum one~\cite{Lig55,Ric56},
which assumes the existence of a steady-state relationship for the traffic,
called the \textit{fundamental diagram}, giving the traffic flow in function of the vehicular density in the road section.
In this paper, we present a reformulation of the $M/G/c/c$ state dependent queuing model~\cite{Smi94} working with flows
rather than speeds, and propose a stochastic traffic model based on this new formulation and inspired from the Godunov
scheme~\cite{God59, Leb96} of the LWR traffic model. We derive a stationary probability distribution of the $M/G/c/c$ state
dependent queuing model on a road section, by considering density-flow fundamental diagrams rather than density-speed ones.
The model assumes a quadratic fundamental diagram which correctly captures the stationary density-flow relationships in both
non-congested and congestion conditions. By this, we consider the traffic demand and supply functions for the section, and
distinguish three road section systems, an open road section, a constrained road section and a closed road section. Finally,
we present a model of two road sections in tandem, derive some performance measures and show on an academic example that the
model captures the dynamics of the road traffic. The model we present here can be used to the analysis of travel times
through road traffic networks (derive probability distributions, reliability indexes, etc),
as in~\cite{ Far08, Far14a, Far14b}, where an algebraic deterministic approach is used. The approach present here remains
also valid with a triangular fundamental diagram (which is the most used), but the formulas will change~\cite{guerrouahane2015}.

The remainder of this article is organized as follows. In Section 2, we first present a review of the existing works.
In this regard, we present in Section 3 a short review of the $M/G/c/c$ state dependent queuing model of
Jain and Smith~\cite{Smi97}. We then propose a new reformulation of this model by working with flows rather than speeds.
We rewrite the stationary probability distribution of the number of cars on the road section.
In Section 4, we present a model of two road sections in tandem where both road sections influence each other.
An implicit system given by an algebraic equation is solved, in order to deal with this mutual dependence.
We derive some performance measures of both road sections (throughput and expected travel time) and compare our
simulation results with those of existing $M/G/c/c$ state dependent queuing models (Kerbache and Smith model). In Section 5, we briefly summarize our findings and future work.

\section{Literature review}

Congestion is a phenomenon that arises on local as well as large areas, whenever traffic demand exceeds traffic supply.
Traffic flow on freeways is a complex process with many interacting components and random perturbations such
as traffic jams, stop-and-go waves, hysteresis phenomena, etc. These perturbations propagate from downstream to upstream sections.
During traffic jams, drivers are slowing down when they observe traffic congestion in the downstream section, causing upstream
propagation of a traffic density perturbation. A link of the network is modeled as a sequence of road sections, for which
fundamental diagrams are given. As the sequence of sections is in series, if any section is not performing optimally,
the whole link will not be operating efficiently. We base here on the LWR first order model~\cite{Lig55, Ric56}, for which numerical
schemes have been performed since decades~\cite{ Dag94, God59, Leb96}.
In the Godunov scheme~\cite{Leb96}, as well as in the cell transmission model~(CTM)~\cite{Dag94}, traffic demand and supply functions
are defined and used.
The demand-supply framework provides a comprehensive foundation for the LWR first order node models. Flow interactions in these models result from
limited inflow capacities of the downstream links. Recently, in~\cite{Tam10} framework has been supplemented with richer features such as
conflicts within the node. In~\cite{Oso11}, a dynamic network loading model that yields queue length distributions was presented. This model is a discretized,
stochastic instance of the kinematic wave model (KWM). In ~\cite{Boel06}, the compositional stochastic model extends the cell transmission
model~\cite{Dag94} by defining demand and supply functions explicitly as random variables, and describing how speed evolves dynamically in each
link of the road network. Several simulation models based on queueing theory have been developed, but few studies have explored the potential of the queueing theory framework to develop analytical traffic models. The development of analytical, differentiable, and computationally tractable probabilistic
traffic models is of wide interest for traffic management. The most common approach is the development of analytical stationary
models~\cite{Oso11}. A review of stationary queueing
models for highway traffic and exact analytical stationary queueing models of unsignalized intersections is given in deferent
literatures~\cite{HeiO1, Van08, Vandael00}. In~\cite{Balsamo,Guerr14, Hei97}, the authors contributed to the study of signalized
intersections and presented a
unifying approach to both signalized and unsignalized intersections. These approaches resort to infinite capacity queues, and thus
fail to account for the occurrence of breakdown and their effects on upstream links.
Modeling and calculating traffic flow breakdown probability remains an important issue when analyzing the stability and reliability of transportation
system ~\cite{Birlon05, Break10}.

Finite capacity queueing theory imposes a finite upper bound on the length of a queue. This allows to account for finite link lengths,
which enables the modeling and analysis of breakdowns. Finite capacity queueing network (FCQN) model are of interest for a variety of
applications such as the study of hospital patient flow, manufacturing networks, software architecture networks, circulation systems,
prison networks~\cite{Smi94}, etc.
The methods in~\cite{Bed12, Smi94, Smi97} resort to finite capacity queueing theory and derive stationary performance measures.
The stationary distribution of finite capacity queueing network exist only for networks with two or three queues with specific topologies~
\cite{Grassman, Langari, Latouche}. An analytic queueing network model that preserves the finite capacity for two queues, under a set
of service rate scenarios proposed in~\cite{Osorio09}.

For more general networks, approximation methods are used to evaluate the stationary distributions, and the General Expansion Method (GEM)
was developed ~\cite{Kerbache1, Kerbache2}. The GEM characterized by an artificial holding node, which is
added and preceded each finite queue in the network in order to register all blocked cars caused by the downstream section. The interested reader may check details in~\cite{Kerbache2}. In~\cite{cruz07}, the authors describe a methodology for approximate analysis of open state
dependent $M/G/c/c$ queueing networks, and a system of two queues in tandem topology was analyzed. The evacuation problem was analyzed using $M/G/c/c$ state dependent queuing networks in~\cite{cruz05} when
an algorithm was proposed to optimize the stairwell case and increase evacuation times towards the upper stories. In~\cite{guerrouahane2015},
the authors present another version of $M/G/c/c$ state dependent queuing model, where a triangular fundamental diagram is considered, which
leads to apply the demand and supply functions in order to analyze a system of two queues in tandem topology.

\section{Review on the M/G/c/c state dependent queuing model}

In this section, we present the $M/G/c/c$ state dependent queuing model~\cite{Smi94,Smi97}.
A link of a road network is modeled with $c$ servers set in parallel, where $c$ is the maximum number of cars that can move on the road (density-like).
The car-speed is assumed to be dependent on the number of cars $n$ on the road, according to a non-increasing density-speed relationship.
Two cases of speed are considered, linear and exponential. In the linear case, we have
\begin{equation}\label{eq-vn}
   v_{n} = v_{1} \left(\frac{c-n+1}{c}\right),
\end{equation}
where $v_{i}$ is the speed corresponding to $i$ cars moving on the road, and $v_1$ is the free flow speed.
In the exponential case, we have
\begin{equation}\label{eq-v}
v_{n} = v_{1} \exp\left[-\left(\frac{n-1}{\beta}\right)^{\gamma}\right],
\end{equation}
in which $\beta$ and $\gamma$ are shape and scale parameters respectively.
$$ \begin{array}{ll}
       \beta & = \frac{a-1}{[\ln(v_{1}/v_{a})]^{1/\gamma}}=\frac{b-1}{[\ln(v_{1}/v_{b})]^{1/\gamma}}. \\
             & \\
       \gamma & = \ln[\frac{\ln(v_{a}/v_{1})}{\ln(v_{b}/v_{1})}]/\ln(\frac{a-1}{b-1}).
   \end{array}$$
The values $a$ and $b$ are arbitrary points, used to adjust the exponential curve~\cite{Smi94,Smi97}.
The arrival process of cars into the link is assumed to be Poisson with rate $\lambda$,
while the service rate of the $c$ servers depend on the number $n$ of cars on
the road. A normalized service rate $f(n)$ is considered, and is taken $f(n) = v_{n}/v_{1} \leq 1$. In the linear case, we have
$f(n) = (c-n+1)/c$. In the exponential case, we have $f(n) = \exp[-((n-1)/\beta)^{\gamma}]$. We notice here that $v_{1}$ is the speed
corresponding to one car in the road (ie. the free speed).

The stationary probability distribution $P_n = Prob~ (N = n)$ of the number of cars $N$ in the $M/G/c/c$ state dependent model have been
developed and shown in~\cite{Smi94} to be stochastic equivalent to $M/M/c/c$ queueing model. Therefore, these probabilities can be written as follows.
\begin{equation}\label{eq1}
 P_n = \frac{(\lambda L/v_{1})^{n}}{\prod_{i=1}^{n}i f(i)} P_0, \qquad
 P_0 = \left(1+\sum_{n=1}^{c}\frac{(\lambda L/v_{1})^{n}}{\prod_{i=1}^{n}i f(i)}\right)^{-1},  \quad n=1,...,c.
\end{equation}
where $L$ is the length of the road section.

From $P_n$, other important performance measures can be easily derived.
\begin{enumerate}
     \item The blocking probability : $ P_{c} = P_{0} (\lambda L/v_{1})^c / \prod_{i=1}^{c} i f(i).$
     \item The throughput : $ \theta = \lambda(1-P_{c}).$
     \item The expected number of cars in the section : $ \bar{N} = \sum_{n=1}^{c}n P_{n}.$
     \item The expected service time : $ W = \bar{N} / \theta.$
\end{enumerate}

In the following section, we rewrite the model presented in this section, by considering a car density-flow relationship
rather than a density-speed one. This reformulation permits us to consider demand and supply functions which model limit conditions
for car-traffic, and by that, to model two road sections in tandem, in Section 4.

\subsection{Reformulation of the M/G/c/c state dependent model}
\label{sec-reformulation}

\begin{figure}[htbp]
    \includegraphics[width=7.5cm]{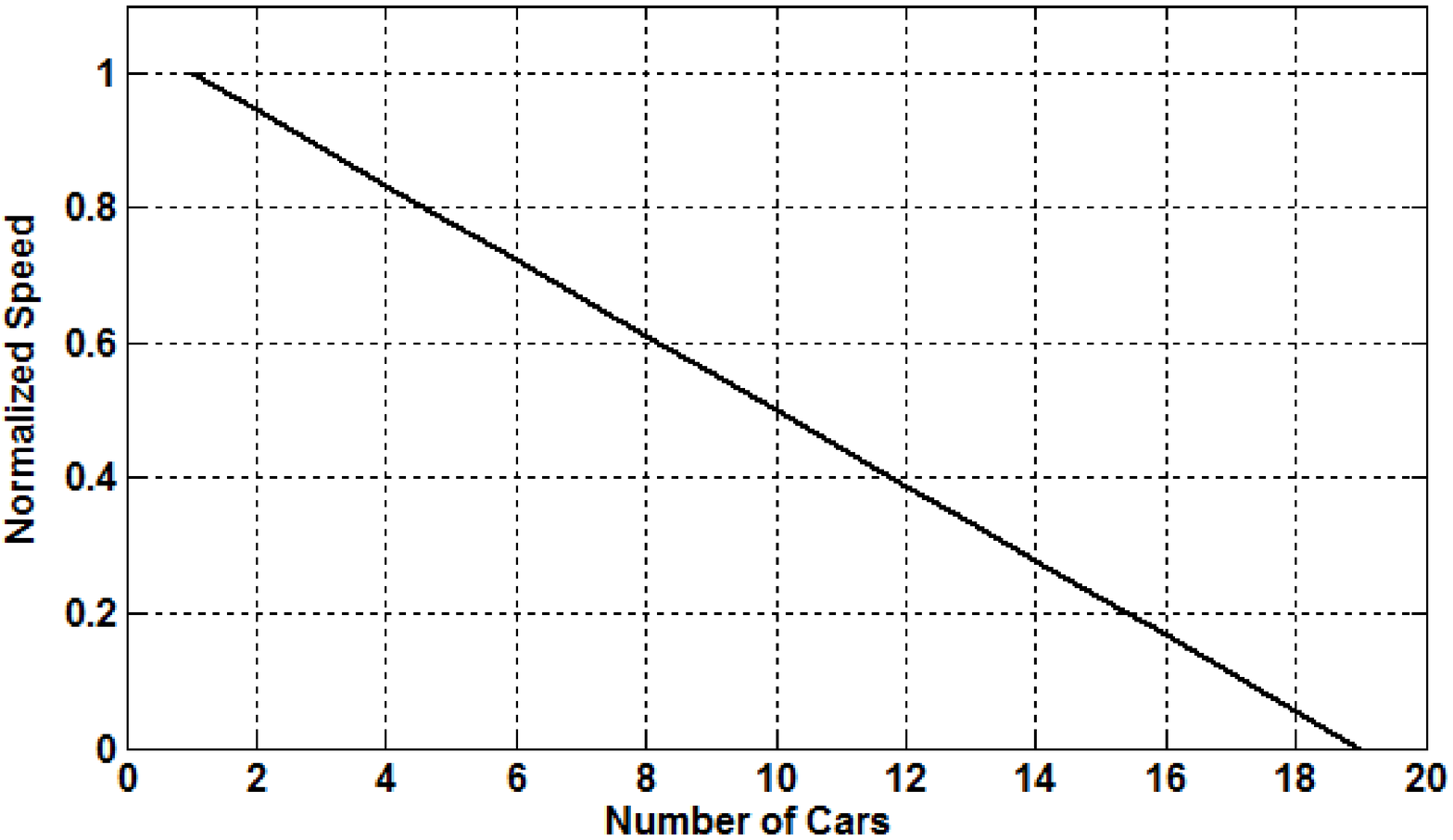}
    \includegraphics[width=7.5cm]{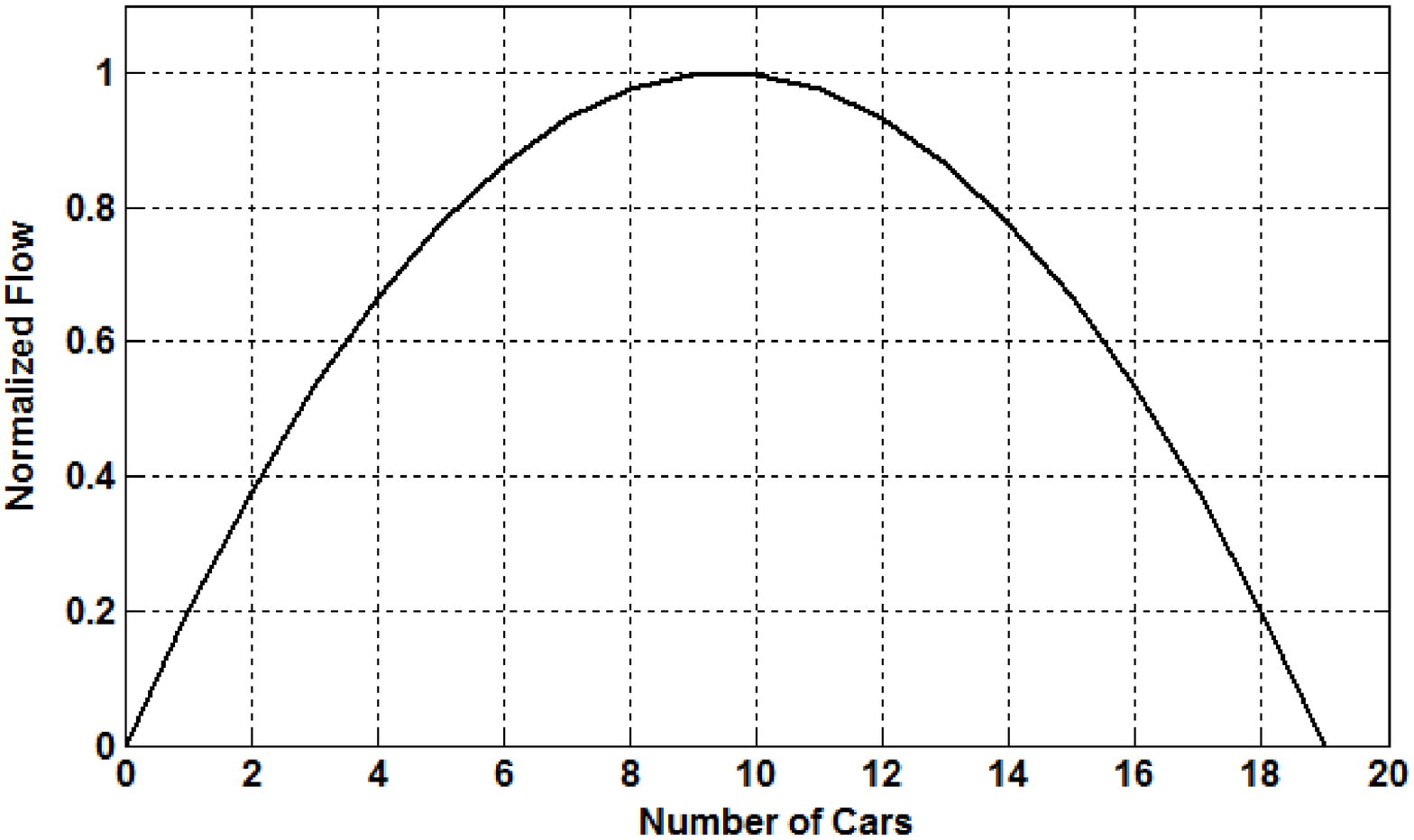}
    \caption{In the left side, the car-density - speed relationship. In the right side, the corresponding car-density - flow relationship.}
    \label{fig-diagram}
\end{figure}

In this section, we first reformulate the $M/G/c/c$ state dependent model of Jain and Smith, by considering car-flows rather than
car-speeds. Indeed, the linear case of that model considers a linear relationship of the speed as a function of the car-density.
This is equivalent to considering a quadratic relationship of the car-flow as a function of the car-density; see Figure~\ref{fig-diagram}.
It is easy to check that this quadratic relationship is written as follows.
\begin{equation}\label{car-flow}
   q_{n} =  q_{\max}\left(1-\left(\frac{c-2n+1}{c+1}\right)^{2}\right).
\end{equation}
Indeed, one can for example verify that parting from~(\ref{car-flow}), we retrieve~(\ref{eq-vn})
\begin{equation} \nonumber
  \frac{v_n}{v_1} := \frac{q_n / \rho_n}{q_1 / \rho_1} = \frac{c-n+1}{c},
\end{equation}
where $\rho_{n}=n/L$.

The stationary probabilities~(\ref{eq1}) are then rewritten as follows.
\begin{equation}\label{eq-pn}
 P_n = \frac{(\lambda / q_{\max})^{n}}{\prod_{i=1}^{n} i g(i)} P_0, \qquad
 P_0 = \left(1+\sum_{n=1}^{c}\frac{(\lambda / q_{\max})^{n}}{\prod_{i=1}^{n} i g(i)}\right)^{-1},  \quad n=1,...,c.
\end{equation}
where
$$g(i) = \frac{q_i/q_{\max}}{i}.$$

The objective of working with flows rather than speeds here is to consider the demand and supply functions which we will use in the following sections,
for the model of two road-sections in tandem.
From~(\ref{eq-vn}) and~(\ref{car-flow}), the relationship between the free speed $v_{1}$ and  the maximum car-flow $q_{\max}$ is given as follows.
\begin{equation}\label{eq3}
q_{\max}=v_{(\frac{c+1}{2})} \rho_{(\frac{c+1}{2})} =\frac{v_{1}}{Lc}\left(\frac{c+1}{2}\right)^{2}.
\end{equation}

The demand $\Delta_{n}$ and the supply $\Sigma_{n}$ functions are given as follows.
\begin{equation}\label{eq-dem}
\Delta_{n} =
\begin{cases}
              q_{n} & \text{if}~~ 0 \leq n \leq \frac{c+1}{2},\\
              q_{\max} & \text{if}~~ \frac{c+1}{2} <n \leq c+1.
\end{cases}
\end{equation}

\begin{equation}\label{eq-sup}
\Sigma_{n} =
\begin{cases}
              q_{\max} & \text{if}~~ 0 \leq n \leq \frac{c+1}{2},\\
              q_{n} & \text{if}~~ \frac{c+1}{2} <n \leq c+1.
\end{cases}
\end{equation}

\begin{figure}[htbp] \label{fig-dem-sup}
    \includegraphics[width=7cm]{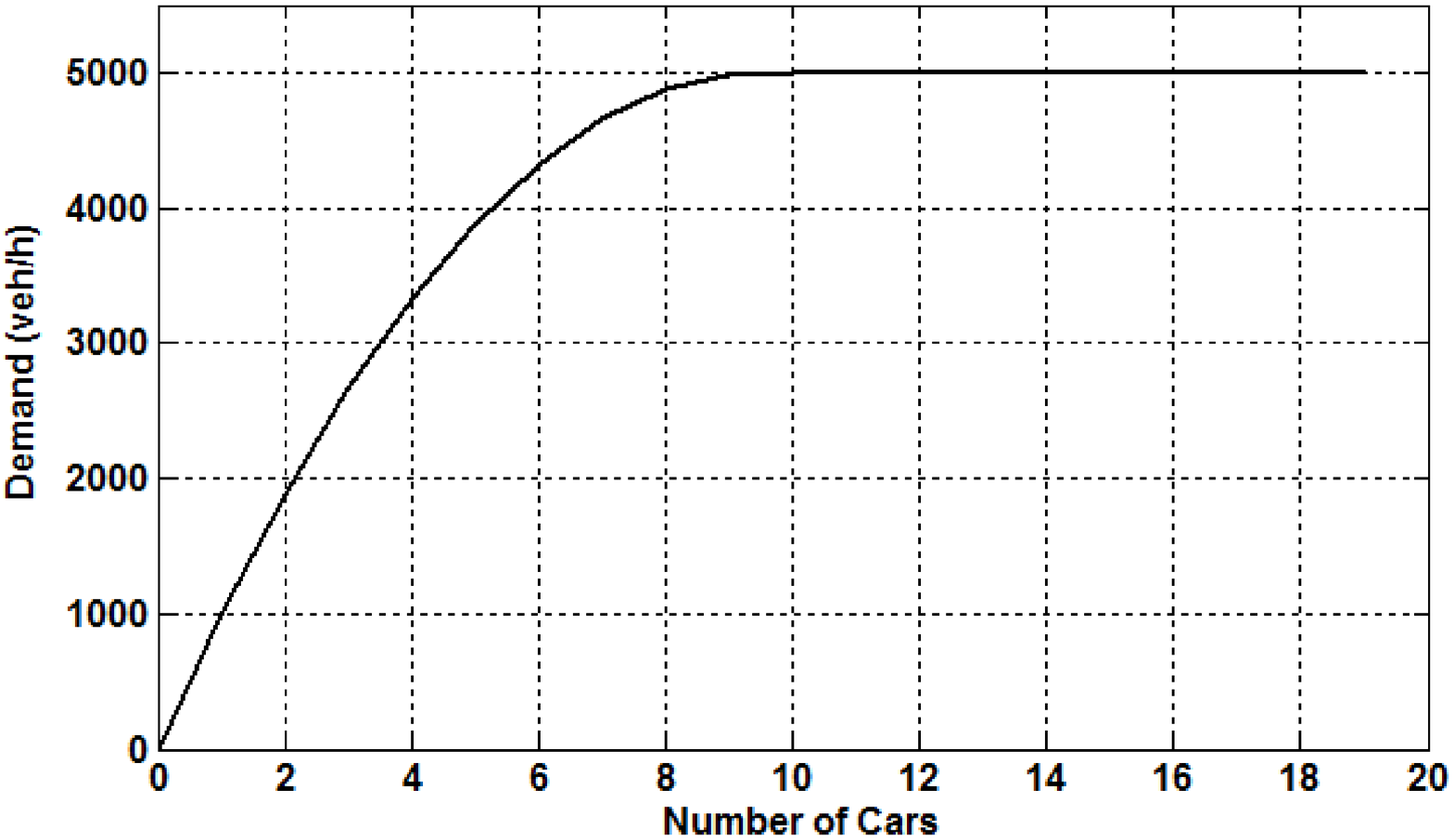}~~~~~
    \includegraphics[width=7cm]{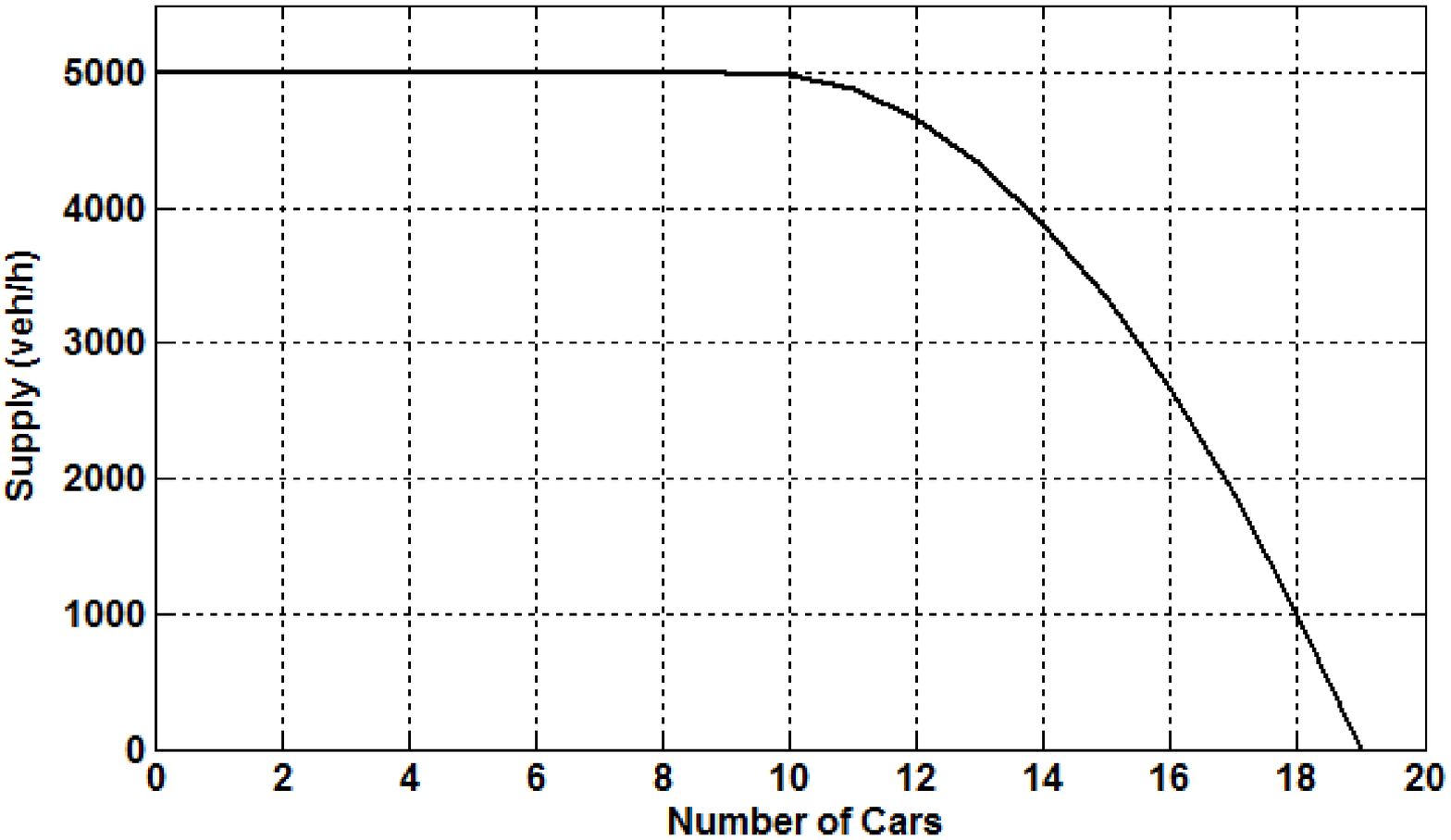}
    \caption{Demand diagram (left). Supply diagram (right).}
    \label{demand-supply}
 \end{figure}

Demand and supply diagrams are shown in Figure~\ref{demand-supply}, where the parameters of Table~1 are considered.

In the following section, we consider three configuration models of a road section: open, constrained and closed road section.
We will explain the role of each model, and tell how these models will be used for the model of two sections in tandem.

\subsection{Open, constrained and closed road sections}

\begin{figure}[htbp]
  \begin{center}
    \includegraphics[width=13cm]{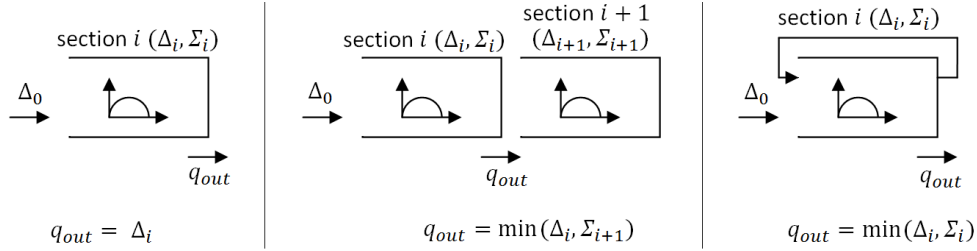}
    \caption{open, constrained and closed road section systems, respectively from left to right sides.}
    \label{fig-occ}
  \end{center}
 \end{figure}

We distinguish here three road section systems.
\begin{itemize}
  \item An open road section is a road section where the traffic out-flow from the section is only constrained by the traffic demand of the section,
    given in function of the number of cars in the section. It is equivalent to the traffic system presented above in Section~\ref{sec-reformulation},
    where the traffic fundamental diagram is restrained to only its traffic demand part (ie with infinite traffic supply).
    In such a road section, which we index by $i$, if we denote by $\Delta_i$ and $\Sigma_i$ the traffic demand and supply of the section
    (induced by the car-density of the section) respectively, and $q_{out}$ the outflow from the section, then we have (see Figure~\ref{fig-occ})
    \begin{equation}
        q_{out} = \Delta_i.
    \end{equation}
  \item A constrained road section is a road section out flowing to another downstream section, and where the outflow of the considered road section
      is constrained by the traffic demand of that section, and by the traffic
      supply of the downstream section. If we index the downstream section by the index $i+1$, denote by $\Sigma_{i+1}$
      the traffic supply of the downstream section, and consider the same notations used above, then we have (see Figure~\ref{fig-occ})
    \begin{equation}
        q_{out} = \min (\Delta_i, \Sigma_{i+1}).
    \end{equation}

  \item A closed road section is the road section system presented above in Section~\ref{sec-reformulation}, for which the outflow
    is given by the fundamental traffic diagram with both demand and supply parts, and where both the parts are
    constrained by the car-density of the section. Using the same notation as above, the out-flow from the section
    is given as follows (see Figure~\ref{fig-occ})
    \begin{equation}
        q_{out} = \min (\Delta_i, \Sigma_i).
    \end{equation}
\end{itemize}

Basing on the three systems presented above, we present in the following section, a model of two road sections in tandem.
It is a kind of a constrained section where the downstream section is a closed one.

\section{Two road sections in tandem}

\begin{figure}[htbp]
\begin{center}
    \includegraphics[width=8cm]{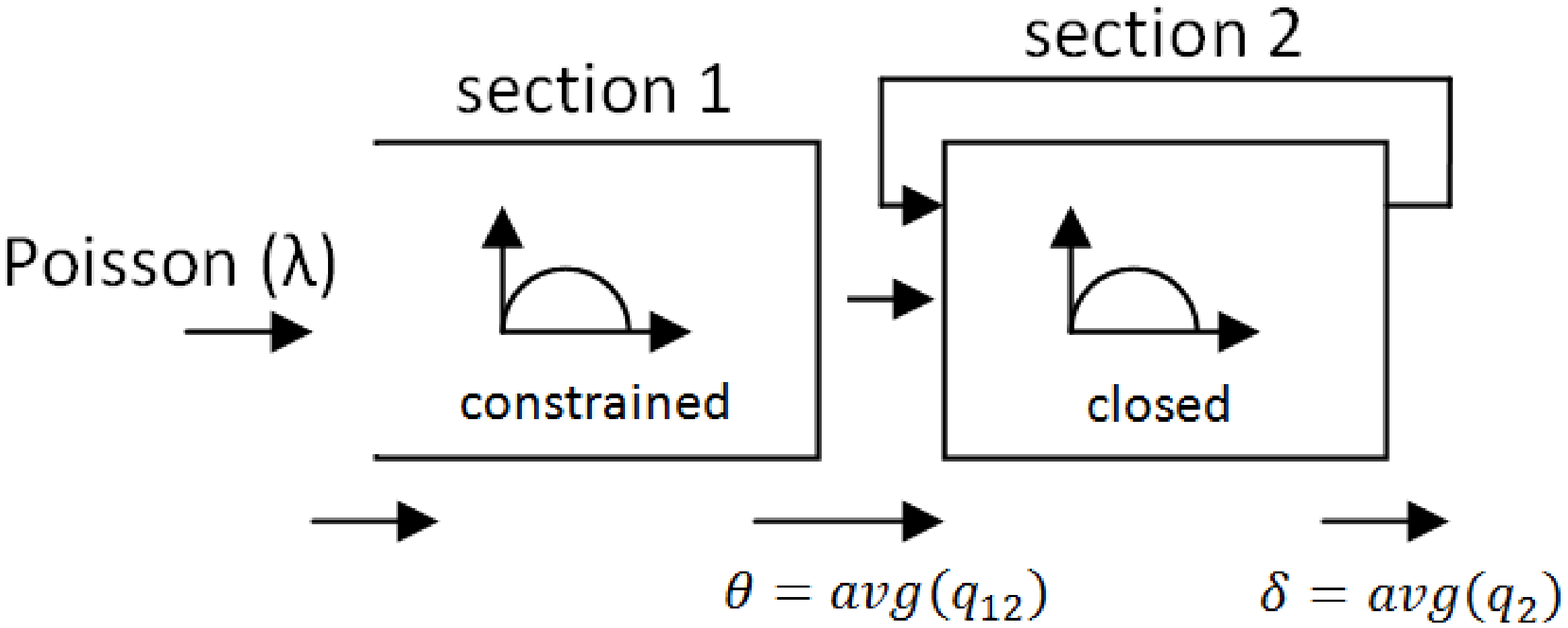}
    \caption{The system of two sections in tandem.}
\end{center}
    \label{tandem}
\end{figure}

We propose in this section a state-dependent $M/G/c/c$ model for two road sections in tandem.
The whole system is a concatenation of two road sections: a $M/G/c/c$ state dependent constrained section (section 1),
with an $M/G/c/c$ state dependent closed section (section 2); see~Figure~\ref{tandem}.
Since the downstream section (section 2) is a closed one, which corresponds to the traffic system presented in Section~\ref{sec-reformulation},
the supply flow of section 2 is stochastic, and is given in function of the probability distribution of the number of cars in section 2.

We assume quadratic flow-density fundamental diagrams for the two sections.
The car-flow outgoing from section 1 and entering to section 2 ($q_{12}(n_1,n_2)$) is given by the minimum between the traffic demand on section 1
and the traffic supply of section 2.
\begin{equation}\label{q12}
    q_{12}(n_1,n_2) = \min(\Delta_1(n_1), \Sigma_2(n_2)).
\end{equation}
The car-flow outgoing from section 2, and also from the whole system is denoted by $q_2$.
The objective here is to determine the stationary probability distributions $P^{(1)}$ and $P^{(2)}$ of the number of cars in section 1 and
section 2 respectively, as well as the stationary probability distribution $P^{(1,2)}$ of the couple $(n_1,n_2)$ of numbers of cars
in sections 1 and 2 respectively.

We notice here that this cannot be calculated directly since both road sections influence each other.
Indeed, $P^{(2)}$ depends on $q_{12}$ which depends on the traffic demand of section 1.
Similarly, $P^{(1)}$ depends on the traffic supply of section 2, since section 1 is constrained by section 2.
We propose below an iterative approach for the calculus of such probabilities.

\subsection{The model}

The procedure we propose here for the determination of $P^{(1)}, P^{(2)}$ and $P^{(1,2)}$ consists in solving an implicit
system given by an algebraic equation, in order to deal with the mutual dependence between sections 1 and 2.
We denote by $\theta$ the average outflow from section 1 (average of $q_{12}$), which is also the average inflow to section 2; and by
$\delta$ the average outflow from section 2 (average of $q_2$), which is the average outflow of the whole system.

We notice that if $\theta$ is known, then the probability distributions $P^{(1)}, P^{(2)}$ and $P^{(1,2)}$ can be obtained in function
of $\theta$. Indeed, section 2 being a $M/G/c/c$ state dependent closed section with an average arrival flow $\theta$,
$P^{(2)}$ can be obtained in function of $\theta$. After that, section 1 being constrained by the supply of section 2 (with known $P^{(2)}$),
$P^{(1,2)}$ and then $P^{(1)}$ can be obtained in function of $P^{(2)}$ and thus in function of $\theta$.
Moreover, from $P^{(1)}$ we can obtain $q_{12}$ in function of $\theta$, and then the expectation of $q_{12}$ gives $\theta$ in function of itself.

From the latter remark, we propose the following iterative algorithm for the calculus of $\theta$ and of the probabilities $P^{(1)}, P^{(2)}$ and $P^{(1,2)}$.

\noindent
\textbf{Step $0$.} Initialize $\theta = \theta_0$. \\
\textbf{Step $k$ (for $k\geq 1$)}.
\begin{itemize}

\item[(1)] \textbf{$P^{(2)}$ is obtained in function of $\theta = \theta_{k-1}$}. \\
      Section 2 is a $M/G/c/c$ state dependent closed section, with an average arrival flow $\theta$, which is unknown.
      According to~(\ref{eq-pn}), the stationary probability distribution of the number $n_2$ of cars on section 2 is given, in function of $\theta$, as follows.
      \begin{equation}\label{eq4}
      \begin{array}{ll}
	  P^{(2)}_{n_2}(\theta) & = \frac{(\theta / q_{\max_2})^{n_2}}{\prod_{i=1}^{n_2}i g_2(i)}  P^{(2)}_{0}(\theta), \qquad n_2=1,...,c_2. \\ \\
	  P^{(2)}_{0}(\theta) & = \left(1+\sum_{n_{2}=1}^{c_{2}}\frac{(\theta / q_{\max_2})^{n_{2}}}{\prod_{i=1}^{n_{2}}i g_2(i)}\right)^{-1}.
      \end{array}
      \end{equation}
      where
      $$g_2(i) = \frac{q_2(i)/q_{\max_2}}{i}.$$

\item[(2)] \textbf{$P^{(1,2)}$ and then $P^{(1)}$ are obtained from $\lambda$ and $P^{(2)}$, and thus in function of $\theta = \theta_{k-1}$}.\\
      Section 1 is constrained by the supply of section 2. As the normalized service rate $g_1(i_1,i_2)$ of section 1
      depends on the number of cars on section 2, which itself depends on $\theta$, the stationary probability distribution $P^{(1)}$ of the number
      of cars in section 1, is also given in function of $\theta$. The normalized service rate $g_1(i_1,i_2)$ of section 1 is given as follows.
      \begin{equation}\label{f12}
		g_1(i_1,i_2) = \frac{q_{12}/q_{\max_1}}{i_1}=  \frac{\min(\Delta_1(i_1), \Sigma_2(i_2))/q_{\max_1}}{i_1}.
      \end{equation}

      In order to write $P^{(1)}$, let us first write
      the stationary probability distribution of the number of cars on section 1, conditioned by the number of cars on section 2, which
      we denote by $P^{(1\mid 2)}$, and which is independent of $\theta$.
      \begin{equation}\label{eqq12}
	\begin{array}{ll}
	  P^{(1\mid 2)}_{(n_{1}\mid n_{2})}(\lambda) & = Prob~(N_1 = n_1 \mid N_2 = n_2) = \frac{(\lambda / q_{\max_1})^{n_1}}{\prod_{i=1}^{n_1} i g_1(i_1, n_2)}  P^{(1\mid 2)}_{(0 \mid n_{2})}(\lambda), \\ \\
	  P^{(1\mid 2)}_{(0 \mid n_{2})}(\lambda) & =  Prob~(N_1 = 0 \mid N_2 = n_2) = \left(1+\sum_{n_1 = 1}^{c_1} \frac{(\lambda / q_{\max_1})^{n_1}}{\prod_{i=1}^{n_1} i g_1(i_1, n_2) }\right)^{-1}.
	\end{array}
      \end{equation}

      Thus, the stationary probability distribution $P^{(1)}$ of the number of cars in section 1, which depends on $\lambda$ and $\theta$, is given as follows.
      \begin{equation}\label{eqq5}
	\begin{array}{ll}
	  P^{(1)}_{n_{1}}(\lambda,\theta) & = \sum_{n_{2}=1}^{c_{2}}P^{(1\mid 2)}_{(n_{1}\mid n_{2})}(\lambda) P^{(2)}_{n_2} (\theta), \\ \\
	  P^{(1)}_{0} (\lambda, \theta) & = \sum_{n_{2}=1}^{c_{2}} P^{(1\mid 2)}_{(0 \mid n_{2})}(\lambda) P^{(2)}_{n_{2}} (\theta).
	\end{array}
      \end{equation}

\item[(3)] \textbf{$q_{12}$ is obtained from $P^{(1)}$, and in function of $\theta = \theta_k$, and finally  $\theta_{k+1} = \mathbb E(q_{12}(\theta_k))$}.\\
      The average outflow from section 1, $\theta$, is given in function of $P^{(1)}$ as follows.
      \begin{equation}\label{theta_old}
	\theta_{k} = \lambda \left( 1 - P^{(1)}_{c_1}(\lambda, \theta_{k-1})\right).
      \end{equation}

\end{itemize}

\noindent
The procedure consists then in solving the following implicit equation  on the scalar $\theta$.
\begin{equation}\label{theta}
   \theta = \lambda \left( 1 - P^{(1)}_{c_1}(\lambda, \theta)\right).
\end{equation}
The solution of that equation, if it exists, and if it is unique, gives the value of the asymptotic average flow passing from
section 1 to section 2, and by that, the stationary probability distributions $P^{(2)}, P^{(1\mid 2)}$ and $P^{(1)}$.
The stationary probability distribution $P^{(1,2)}$ can be easily deduced from $P^{(1\mid 2)}$ and $P^{(2)}$.
We notice here that equation~(\ref{theta}) is a fixed-point-like equation.

In the following section, we give some results on the fixed-point equation solving.

\subsection{Fixed-point equation solving}

We use the following notations.
\begin{itemize}
 \item $h(\theta) = \lambda \left( 1 - P^{(1)}_{c_1}(\lambda, \theta)\right)$,
 \item $e(\theta) = h(\theta) - \theta$.
\end{itemize}
Equation~(\ref{theta}) is written $\theta = h(\theta)$, and it is equivalent to $e(\theta) = 0$.

Below, we give some results on the resolution of equation~(\ref{theta}).

\begin{theorem}\label{thm1}
   A solution for equation~(\ref{theta}) exists and is unique.
\end{theorem}
\proof See \ref{appA}.

The following result gives a condition on the couple $(\lambda,\theta)$, under which the sequence $\theta_k, k\in\mathbb N$
defined by fixed point iteration~(\ref{theta_old}) (i.e. $\theta_{k} = h(\theta_{k-1}))$ converges for every initial value $\theta_0\in[0, \lambda]$.

\begin{theorem}\label{thm2}
   If $\exists \varepsilon > 0$ such that $\varepsilon < \theta / \lambda$, and if
   $$\sum_{n_2=1}^{c_2} P^{(1\mid 2)}_{c_1\mid n_2} (\lambda) P^{(2)}_{n_2}(\theta)\left( n_2 - \sum_{k_2=1}^{c_2} k_2 P^{(2)}_{k_2}(\theta)\right) \leq \varepsilon,$$
   then the fixed point iteration converges to the unique fixed point of the fixed point equation.
\end{theorem}
\proof See \ref{appB}.

\begin{figure}[htbp] \label{func h}
  \begin{center}
     \includegraphics[width=9cm]{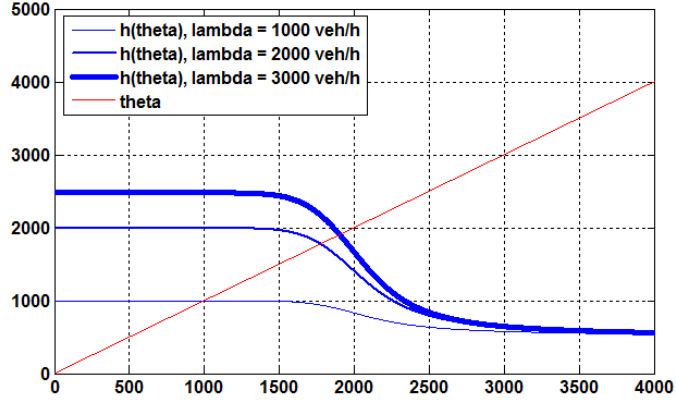}
     \caption{On the x-axis: $\theta$. On the y-axis: in blue the function $h(\theta)$ for different values of $\lambda$: 1000, 2000 and 3000; in red: $\theta$.}
      \label{func-h}
   \end{center}
\end{figure}

\begin{figure}[htbp]
    \includegraphics[width=7.5cm]{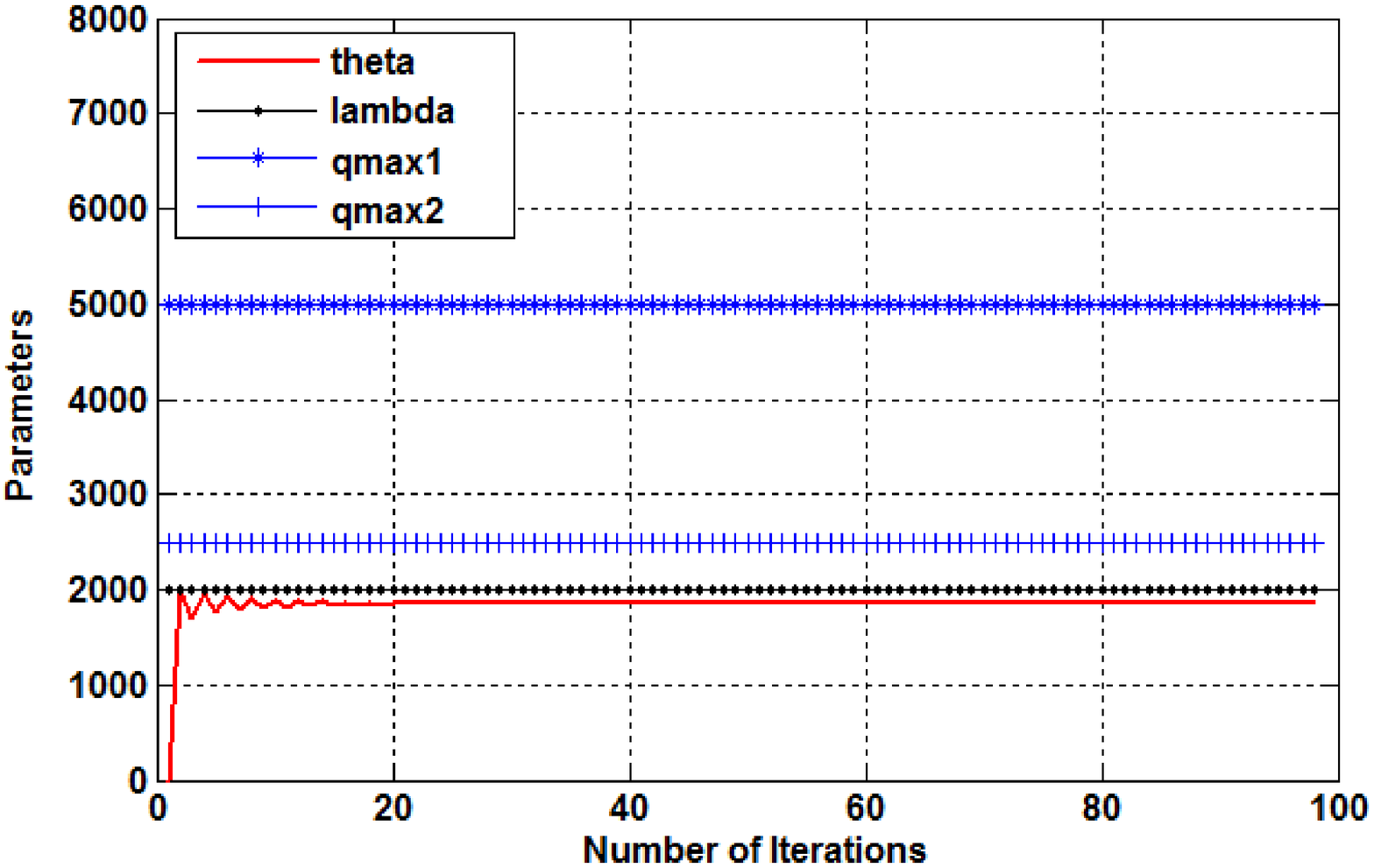}~~~~~~
    \includegraphics[width=7.5cm]{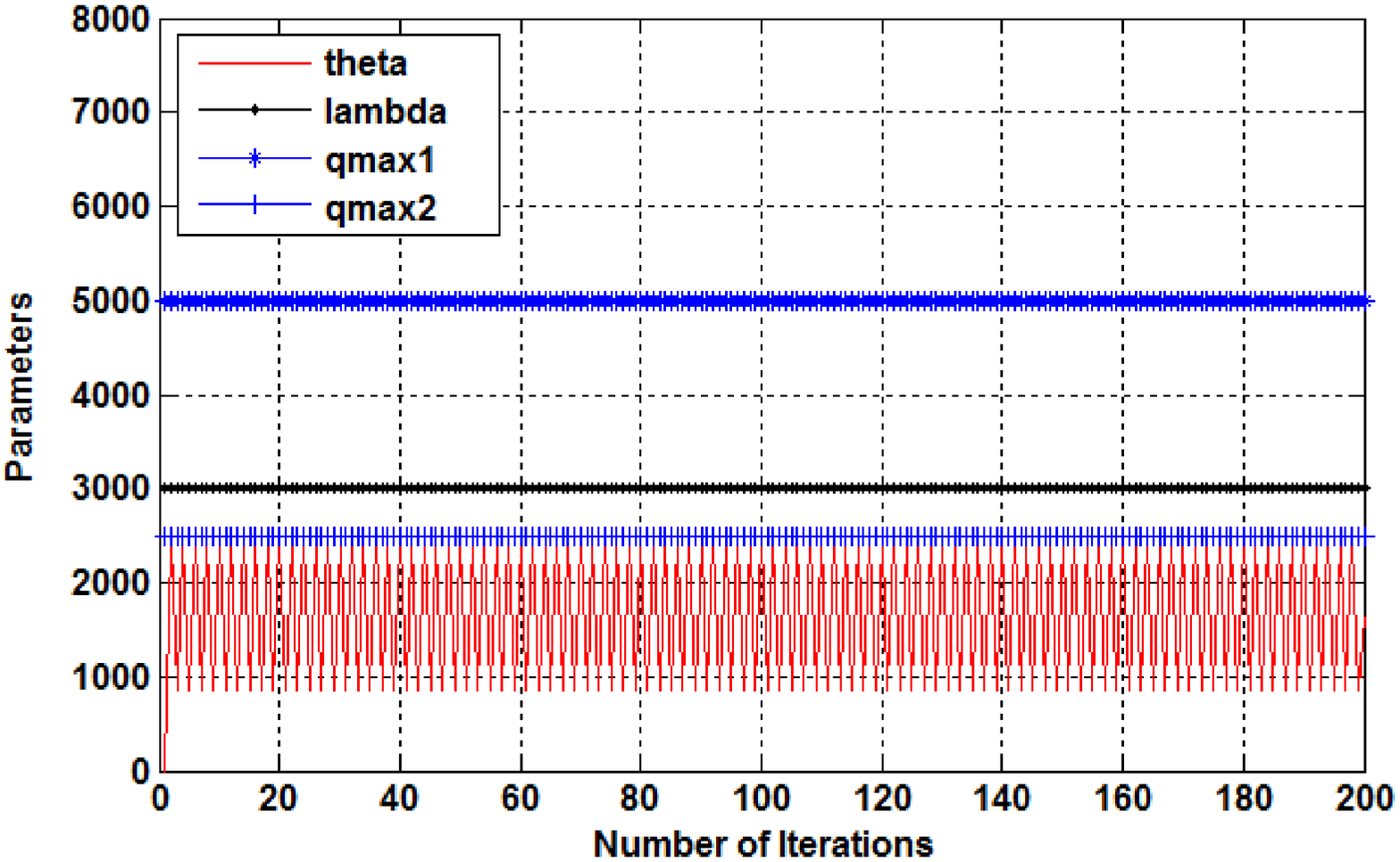}
    \caption{Fixed point iteration. On the left side: $\lambda=2000$ veh/h and the iteration converges to the fixed point. On the right side: $\lambda=3000$ veh/h and
       the iteration oscillates between two values $h(0) = \lambda$ and $h(\lambda)$.}
    \label{cvg}
 \end{figure}

We do not have other explicit formulas than the condition of Theorem~\ref{thm2}.
We know that for low values of $\lambda$, we have $-1 < d\;h(\theta)/d\theta \leq 0, \forall \theta\in [0, \lambda]$,
and the iteration of the fixed point
equation converges, independently of the initial value of $\theta$, to the unique fixed point.
However, we also know that for high values of $\lambda$, there exists an interval $(\theta_1, \theta_2)$ where
$d\;h(\theta)/d\theta < - 1$. In such cases, the fixed point is unstable. Nevertheless, the sequence $\theta_k, k\in\mathbb N$
oscillates between two adherence values $h(0) = \lambda$ and $h(\lambda)$, giving $(\lambda+h(\lambda))/2$ as the asymptotic
average value of $\theta$, see Figures~\ref{func-h} and~\ref{cvg}.

In Figure~\ref{func-h}, we give in blue color the function $h(\theta)$ for different values of $\lambda$.
The red line is the first bisector. The fixed point verifying $h(\theta)=\theta$ is then the intersection of the blue curve with the red one.
We see in this figure that the derivative of $h(\theta)$ with respect to $\theta$ on the fixed point is decreasing with $\lambda$, starting from zero.
Indeed, for low values of $\lambda$, we have $d h(\theta)/\theta = 0$. As $\lambda$ increases, $d h(\theta)/\theta$ decreases.
As explained above, while $d h(\theta)/\theta > -1$, we know that the fixed point is stable. However, once $d h(\theta)/\theta > -1$, the fixed point
is unstable.\\

Figure~\ref{cvg} illustrates the two cases of stable and unstable fixed points, depending on the value of $\lambda$.
In the first case (left side of the figure), we have $\lambda=2000$ veh/h; see also Figure~\ref{func-h}.
In this case, the fixed point is stable, and the iteration converges to it.
In the second case (right side of Figure~\ref{cvg}), we have $\lambda=3000$ veh/h; see also Figure~\ref{func-h}.
In this case, the fixed point is unstable, and the iteration oscillates from the two values $h(0) = \lambda$ and $h(\lambda)$.

In the following section, we explain how to compute the stationary probability distributions in function of the
average flow $\theta$, calculated by the fixed point iteration.

\subsection{The stationary probability distributions}
\label{subsec}

Once $\theta$ is obtained, the whole stationary regime of the system of two road sections in tandem is determined.
Indeed, the stationary probability distributions $P^{(1)}, P^{(2)}$ and $P^{(1|2)}$ are given by~(\ref{eqq5}), (\ref{eq4})
and~(\ref{eqq12}) respectively.
The stationary joint probability distribution $P^{(1,2)}$ of the couple $(n_1, n_2)$ of the number of cars in sections~1
and~2 respectively, is then easily deduced.
$$P^{(1,2)}_{(n_1,n_2)} = P^{(1|2)}_{n_1|n_2} P^{(2)}_{n_2}, \quad \forall 0\leq n_1\leq c_1, 0\leq n_2\leq c_2.$$

In the following, we illustrate those probability distributions for given parameters of the system.
Table~1 gives the fixed parameters for the illustrations, where $L$, $v_1$, $\rho_j$ and $q_{\max}$
denote respectively section length, free speed, jam-density and maximum car-flow.
\begin{table}[htbp]
  \begin{center}
  \caption{Parameters for sections 1 and 2.} ~~\\~~
  \begin{tabular}{c|c|c|c|c}
     Section $i$ & $L$ (km) & $v_{1}$ (km/h) & $\rho_j$ (veh/km) & $q_{\max}$ (veh/h) \\ 
     \hline
     1 & 0.1 & 100  & 180 & 5000 \\ 
     \hline
     2 & 0.1 & 50  & 180 & 2500 \\ 
  \end{tabular}
  \end{center}
\end{table}

Let us notice that $\rho_{cr} = \rho_j /2$ (quadratic fundamental diagram), and $c = L \rho_j$.\\
$\rho_{cr}$ and $c$ denote respectively critical car-density and road section capacity
in term of maximum number of cars.

Figure~\ref{fig-compar1} compares the stationary distribution probability of the two following cases.
\begin{itemize}
  \item One closed road section (Jain and Smith model), with parameters of section 2 in Table 1. (Blue color in Figure~\ref{fig-compar1}).
  \item One constrained road section, with parameters of section 1 in Table 1. The section is assumed to be
    constrained by a closed road section whose parameters are those of section 2 in Table 1.
    (Red color in Figure~\ref{fig-compar1}).
\end{itemize}
The average car inflow $\lambda$ is varied from one illustration to another in Figure~\ref{fig-compar1}.

\begin{figure}[htbp] \label{fig-compar1}
  \includegraphics[width=10cm]{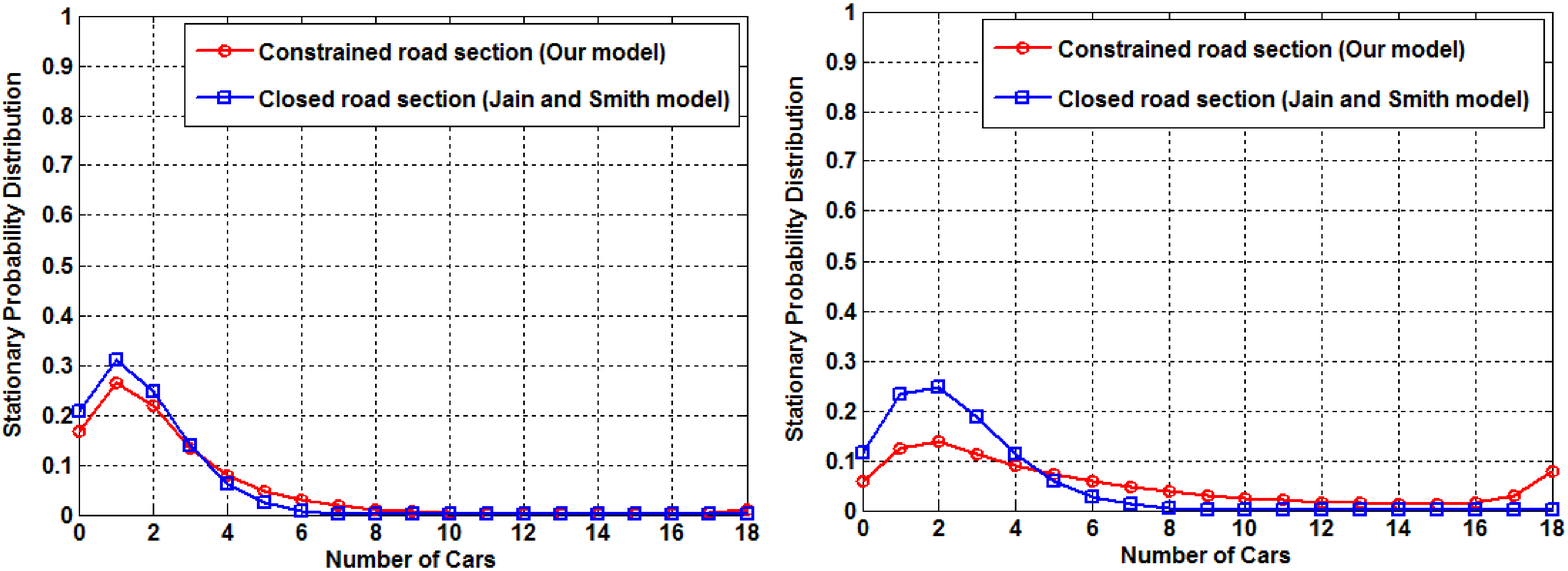}
  \includegraphics[width=5cm]{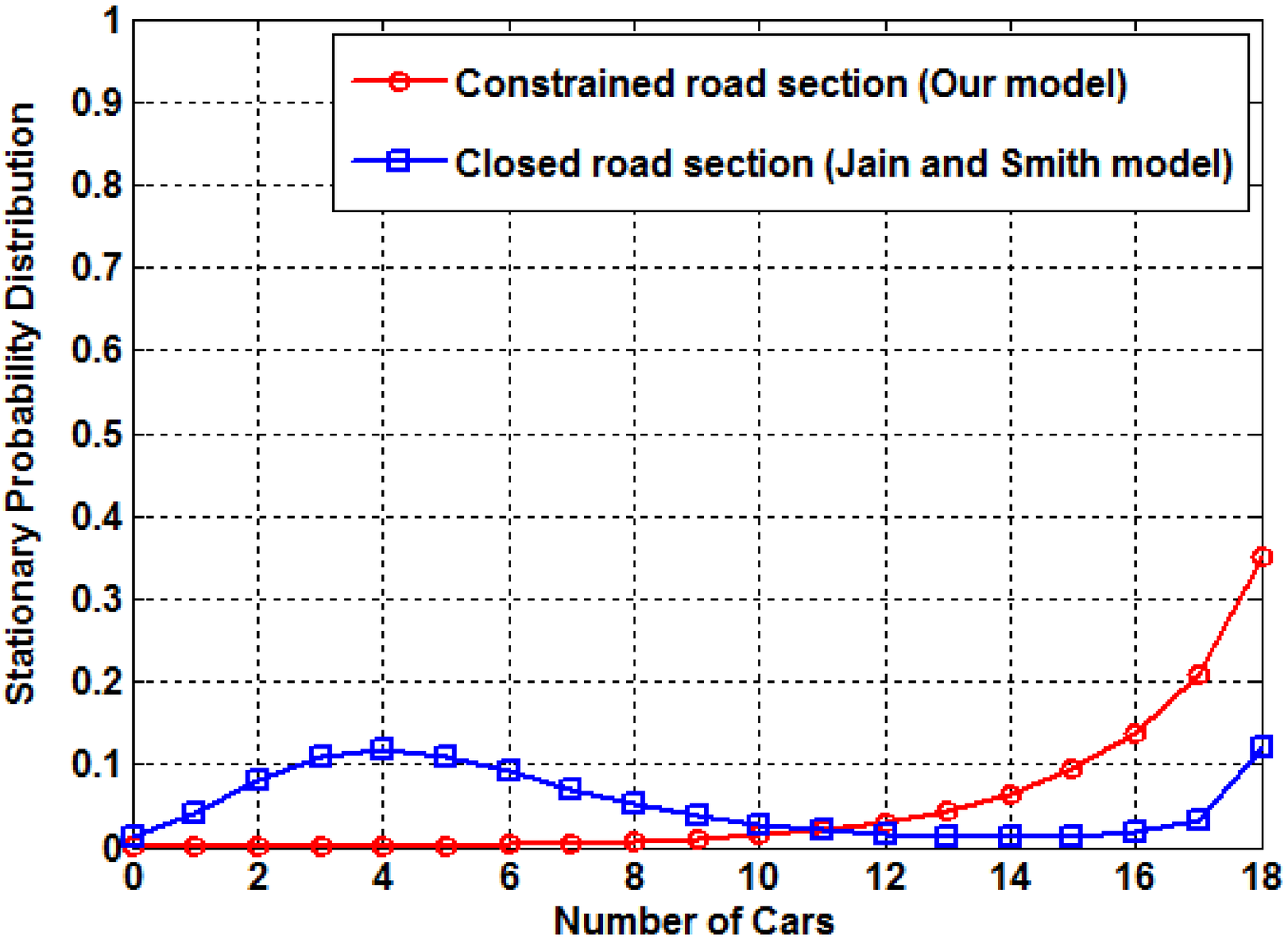}
  \caption{Stationary probability distributions. $\lambda =$ $1000$ veh/h, $2000$ veh/h and $3000$ veh/h respectively from left to right sides.}
  \label{fig-compar1}
\end{figure}

We notice here that the constrained road has a maximum flow capacity of $5000$ veh/h, but it is constrained by a closed
road with a maximum flow capacity of $2500$ veh/h. We remark from Figure~\ref{fig-compar1}, that the constrained
section is more likely to be congested than the closed one, even though a priori they are both constrained
by a maximum flow capacity of $2500$ veh/h.

\begin{figure}[htbp]
     \includegraphics[width=7.5cm]{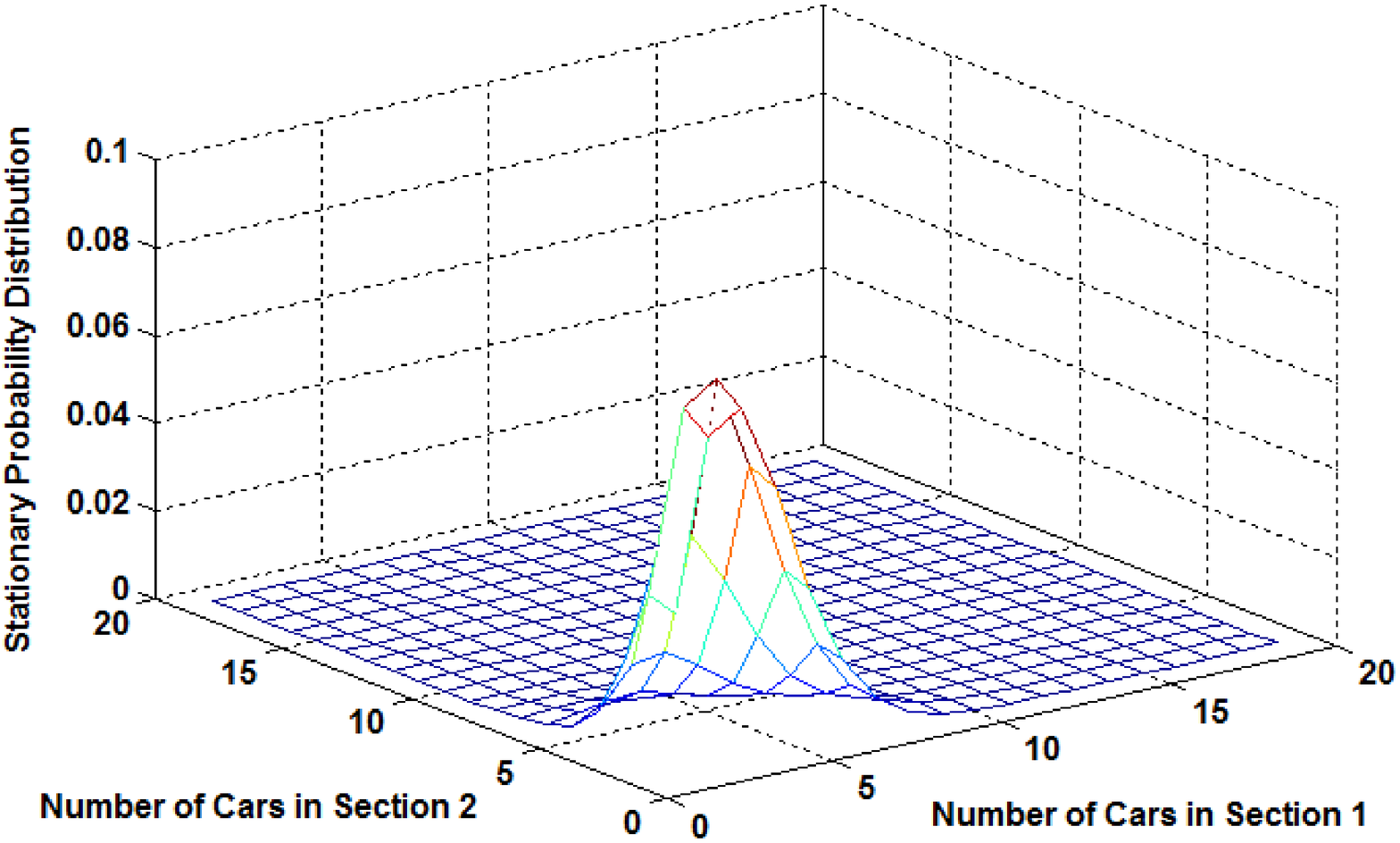}
     \includegraphics[width=7.5cm]{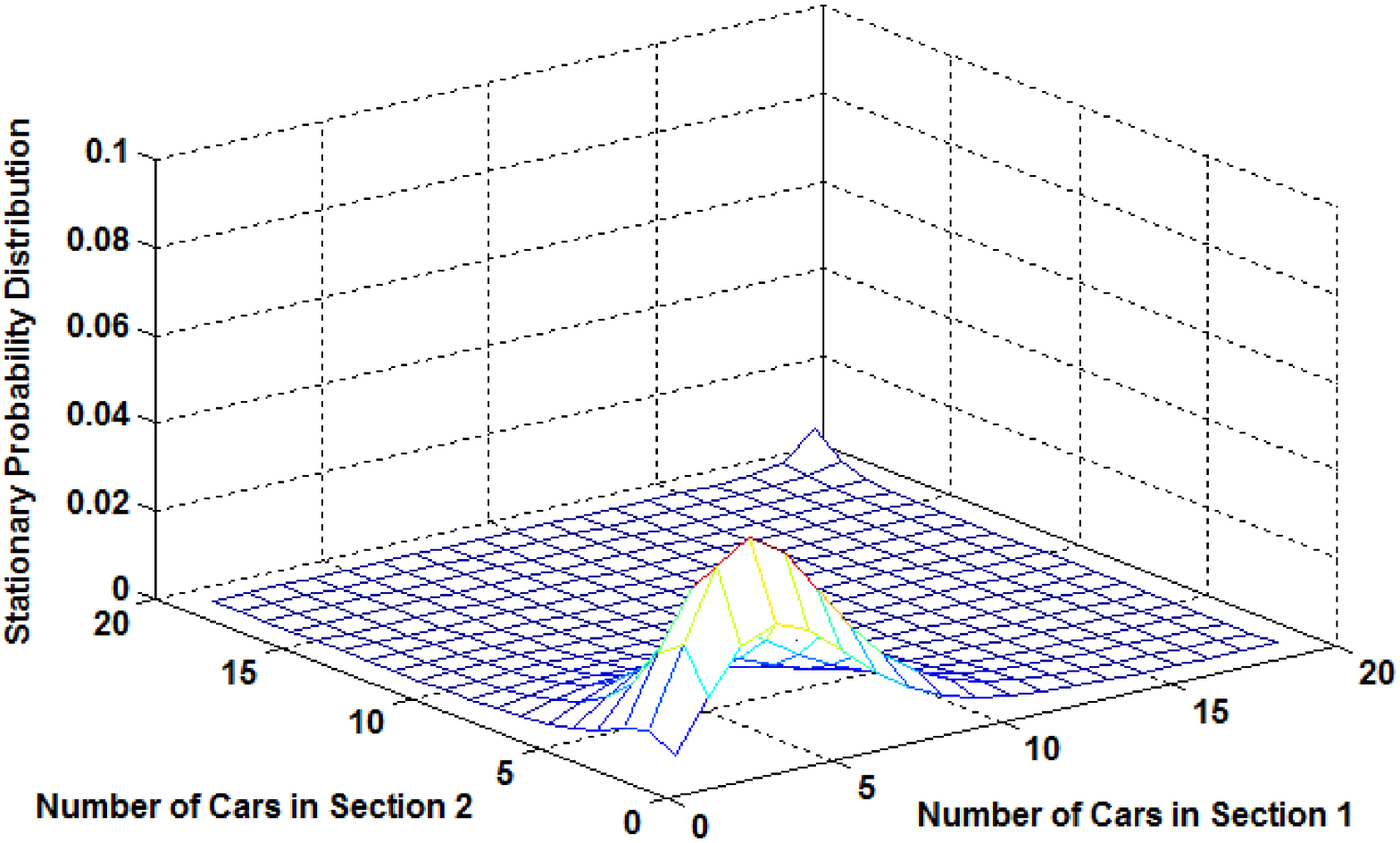}\\
     \includegraphics[width=7.5cm]{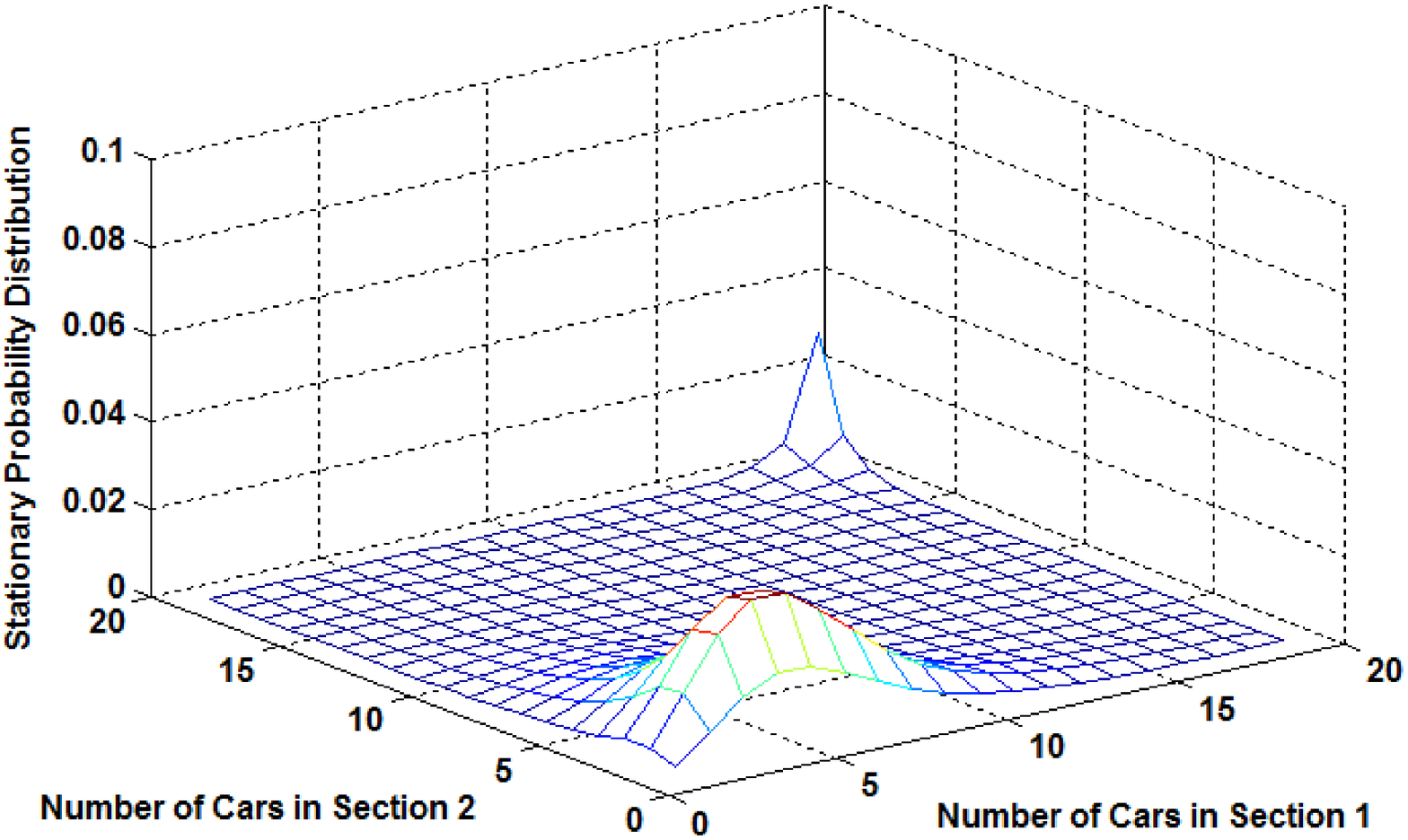}
     \includegraphics[width=7.5cm]{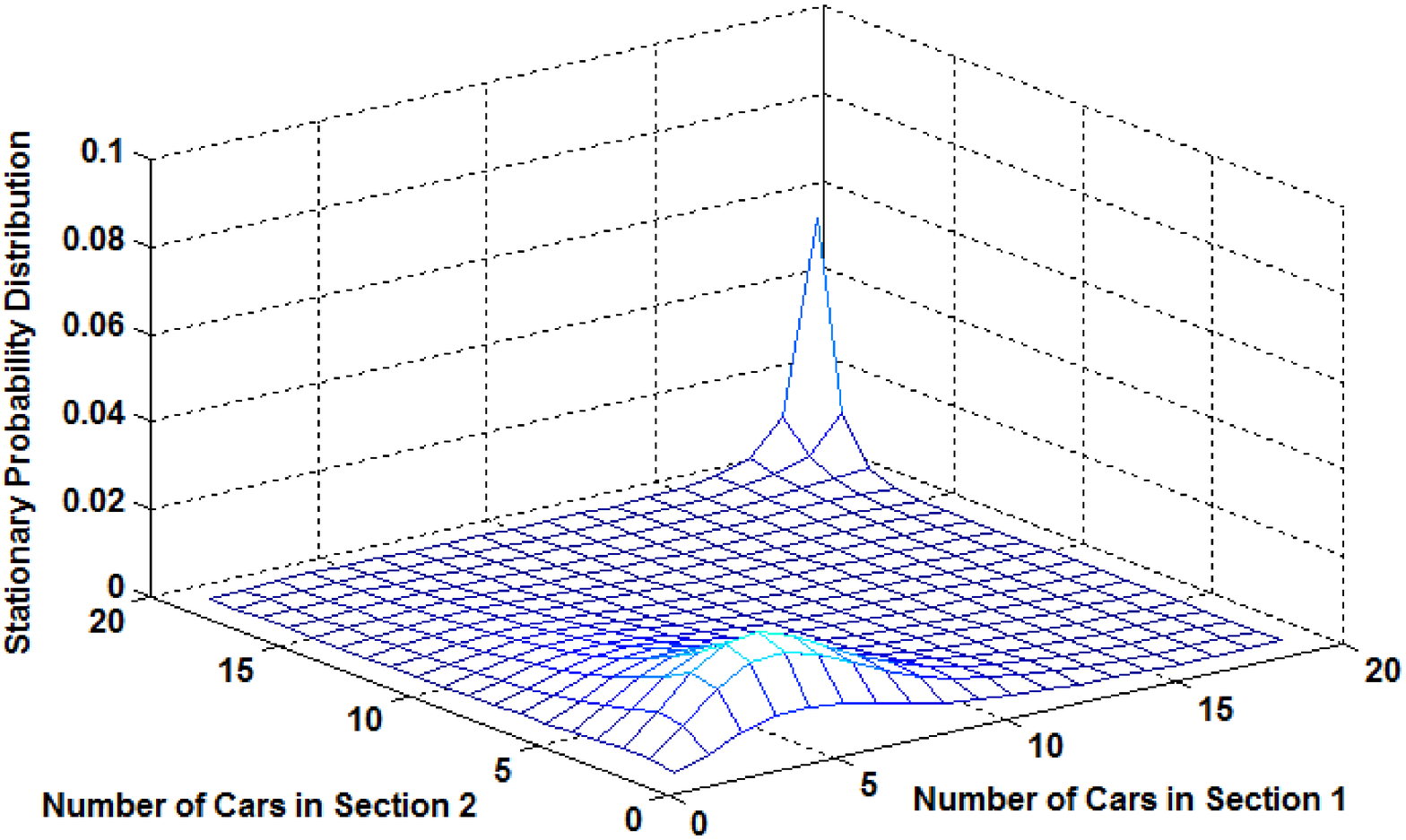}
     \caption{Stationary probability distribution of the number of cars in two road sections in tandem.
        The arrival rate $\lambda =$ 1000 veh/h, 2000 veh/h, 2500 veh/h and 3000 veh/h.}
     \label{distribution-two-section}
\end{figure}

Figure~\ref{distribution-two-section} gives the stationary probability distributions $P^{(1,2)}$ of the couple $(n_1,n_2)$ of numbers of cars
in sections 1 and 2 respectively. The average arrival rate $\lambda$ is varied from one illustration to another.
Parameters of Table~1 are used. $\lambda$ takes the values of $1000$ veh/h, $2000$ veh/h, $2500$ veh/h and $3000$ veh/h.\\

In Figure~\ref{fig-throughput}, we give the throughputs of sections 1 and 2 in the case where the two sections are set in tandem.
The throughput $\theta$ of section 1 is given by~(\ref{theta}), and it is the solution of the fixed point iteration.
The throughput $\delta$ (see~\ref{tandem}) is computed by the same way (Little's law), once $\theta$ is determined.
\begin{equation}\label{delta} \nonumber
   \delta= \theta \left( 1 - P^{(2)}_{c_2}(\theta)\right).
\end{equation}

Figure~\ref{fig-throughput} displays for an increasing arrival rate, $\lambda$, the throughput through sections 1 and 2 of the whole system of
two road sections in tandem. The parameters of each road section are those given in Table~1. We
remark from Figure~\ref{fig-throughput}, that throughput through the sections 1 and 2 is equal to the arrival rate $\lambda$ up to roughly $2200$ veh/h. This comes directly from the fact that the capacities are not restrictive. However, the trend changes at about $2300$ veh/h (well below $2500$ veh/h), and the throughput does not increase anymore.

\begin{figure}[htbp]
  \begin{center}
    \includegraphics[width=8cm]{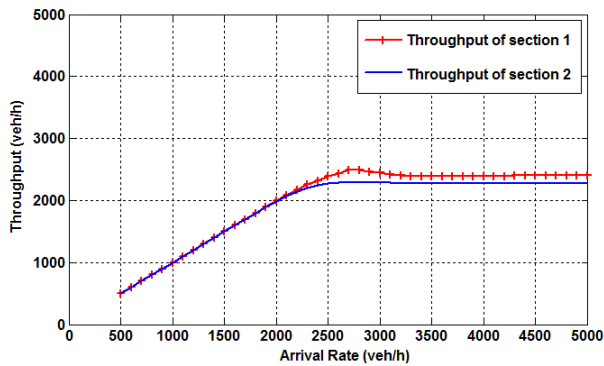} \label{fig-throughput}
    \caption{Throughputs $\theta$ in red color and $\delta$ in blue color in function of the average arrival rate $\lambda$. The two road sections are set in tandem.}
    \label{fig-throughput}
  \end{center}
\end{figure}

\subsection{Performance measures of two sections in tandem topology}

We consider the system of two sections in tandem topology presented in Figure~\ref{tandem}, with parameters of Table~1. The bottleneck in section 2 is most difficult to be analyzed because of the mutual dependence between the two sections set in tandem. In the next, we compare the throughput and the expected service time of our method with the method of Kerbache and Smith~\cite{Kerbache2}.

\begin{figure}[htbp]
  \includegraphics[width=7.5cm]{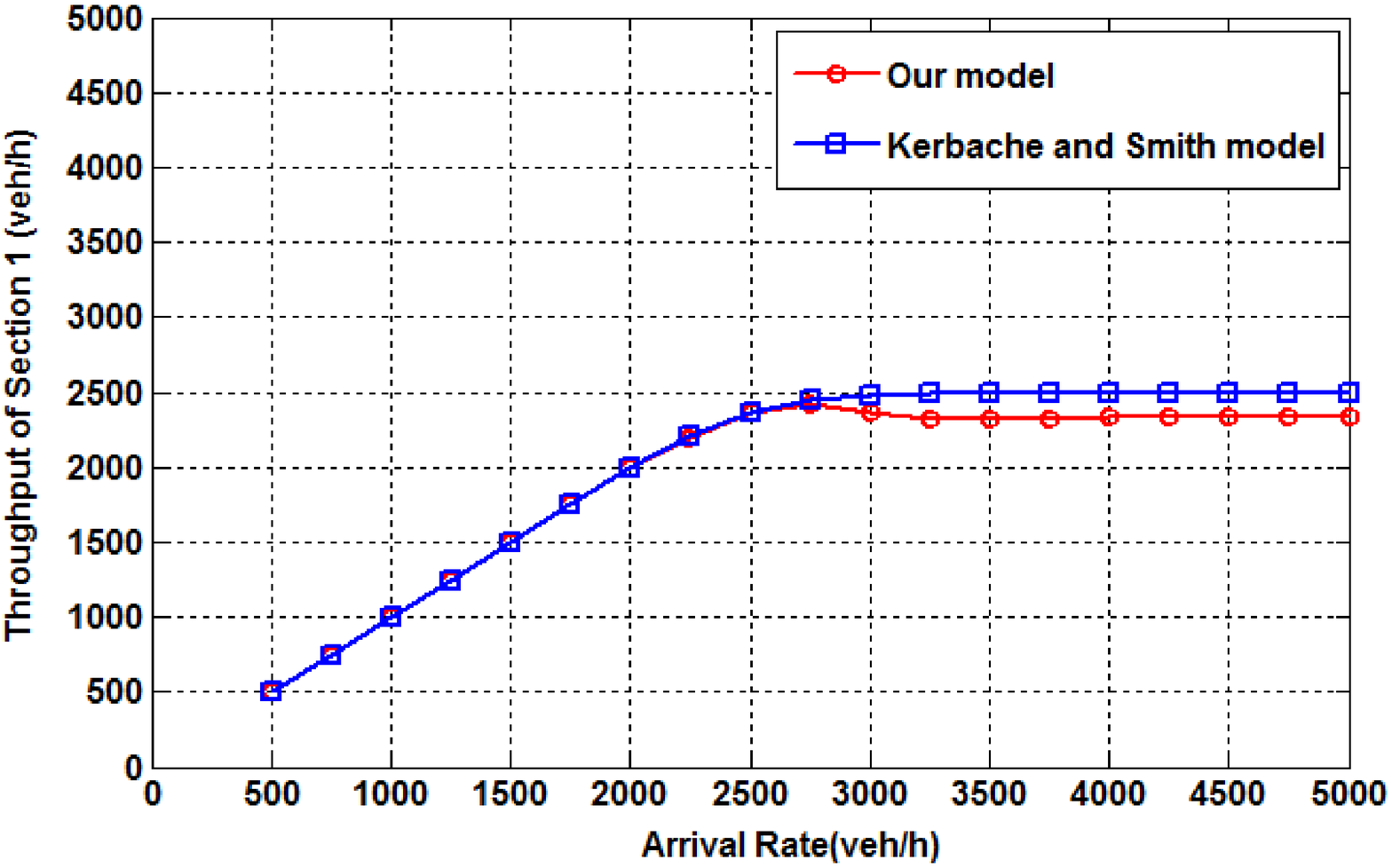}~~~
  \includegraphics[width=7.6cm]{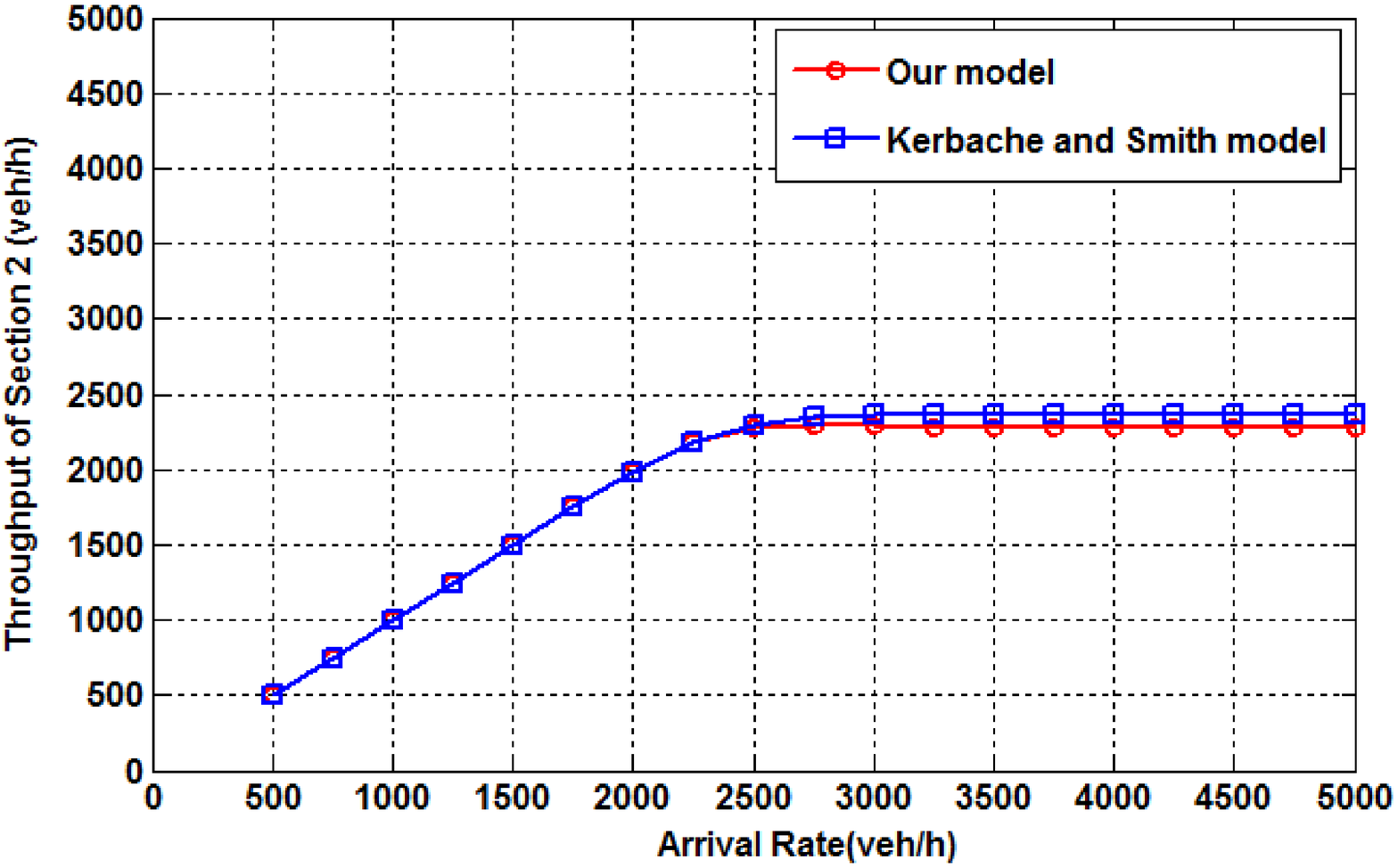}
  \caption{Comparison of the throughputs. In left, section 1 ($\theta$). In right, section 2 ($\sigma$). In red color, our model of two sections in tandem.
     In blue color, open $M/G/c/c$ state dependent queueing networks in tandem}
  \label{fig-compar-throughput}
\end{figure}

Figure~\ref{fig-compar-throughput} compares the throughputs of sections 1 and 2 of our method with the expansion method of Kerbache and Smith
(open $M/G/c/c$ state dependent queueing networks)~\cite{Kerbache2}.\\

We remark from Figure~\ref{fig-compar-throughput}, that when the arrival rate $\lambda$ is low,
our throughputs ($\theta$ of section 1 and $\sigma$ of section 2) are similar to the throughputs of the expansion method of Kerbache and Smith.
In this case, the arrival rate is very light and easily accepted by section 2 without any blocking. When the arrival rate $\lambda$ is large,
our throughputs is low as compared to the throughputs of $M/G/c/c$ state dependent queueing networks, because the constrained section
(section 1 of our model ) is more likely to be congested than the open section (section 1 in the expansion method of Kerbache and Smith).

\begin{figure}[htbp]
  \includegraphics[width=7.5cm]{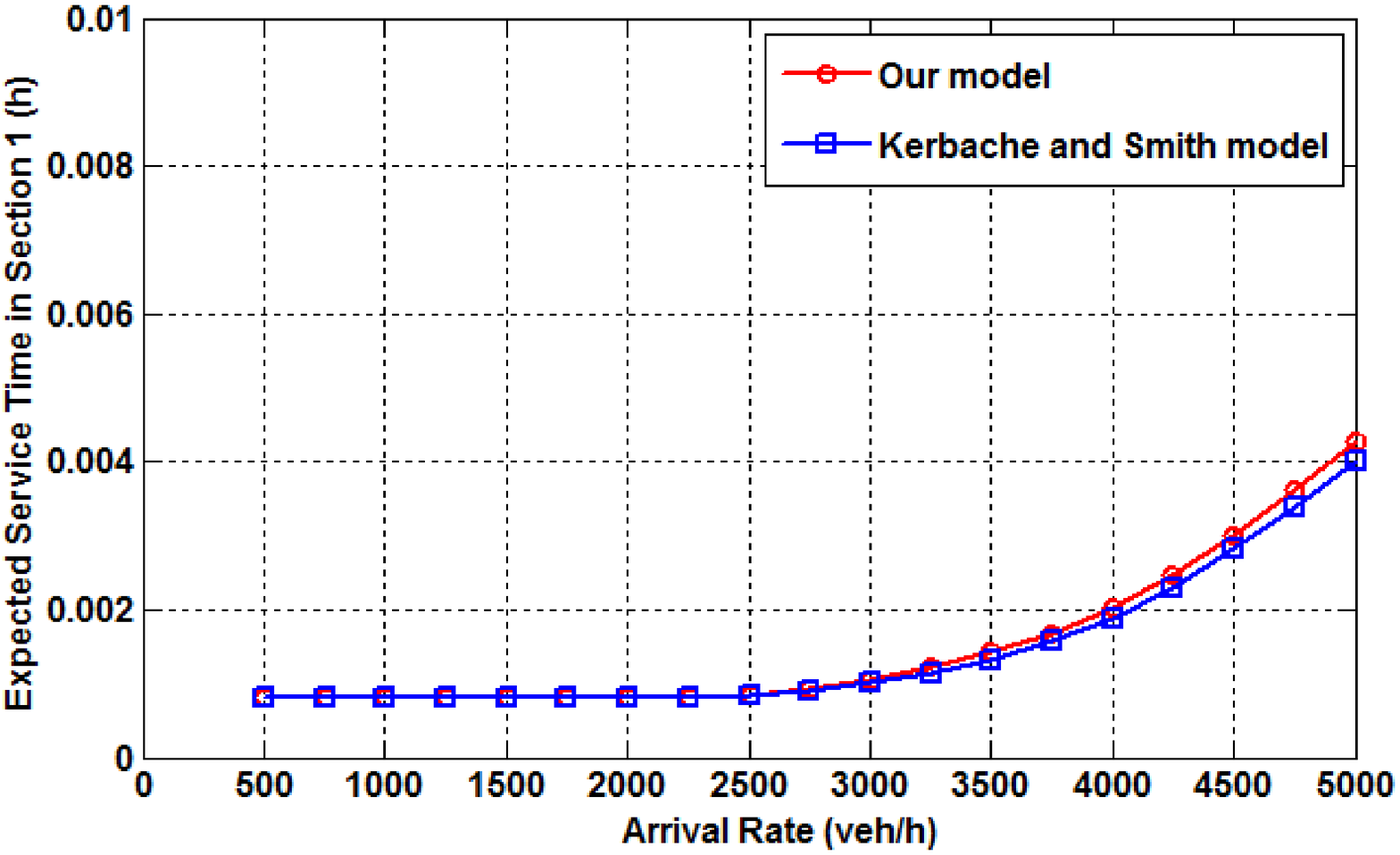}~~~~
  \includegraphics[width=7.6cm]{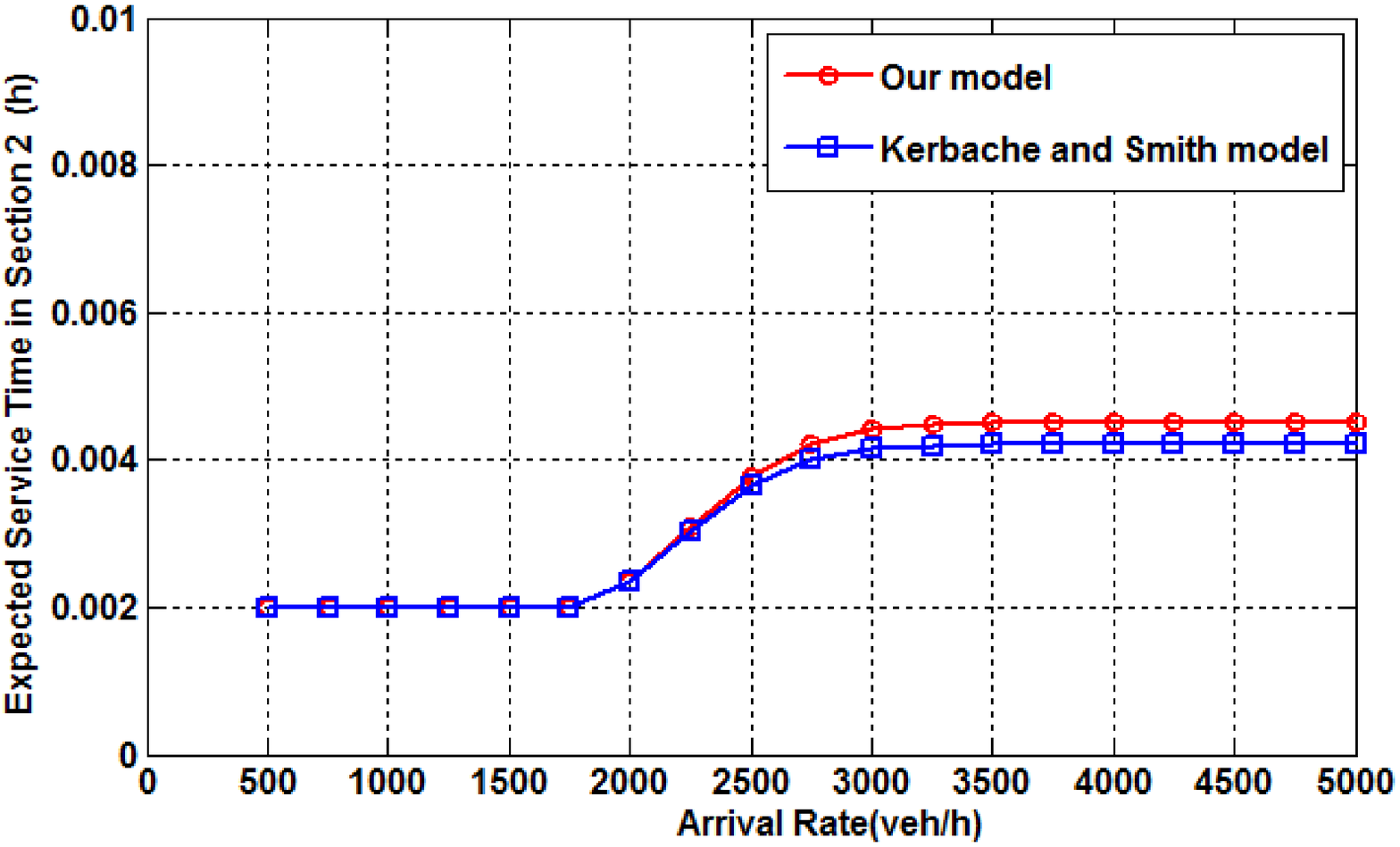}
  \caption{Comparison of the expected service time. In left, section 1. In right, section 2. In red color, our model of two sections in tandem. In blue color, open $M/G/c/c$ state dependent queueing networks in tandem.}
  \label{fig-compar-service-time}
\end{figure}

Using Little's Law, Figure~\ref{fig-compar-service-time} compares the expected service time of sections 1 and 2 of our method with the one of method of Kerbache and Smith~\cite{Kerbache2}.
We remark here that under a light arrival rate $\lambda$, the expected service time almost corresponds to the free service time $(L/v_{1})$.
Under a heavy arrival rate, vehicles slow down and the expected service time increases in value for the two models.
Our expected travel time remains near the upper bound of the expansion method of Kerbache and Smith.

\section{Conclusion and future work}

In this paper, we have presented a queuing model for road traffic that preserves the finite capacity property of the real system.
Based on the $M/G/c/c$ state dependent queuing model, we have proposed a stochastic queuing model for the road traffic
which captures the stationary density-flow relationship in both uncongested and congestion conditions. First results
of the proposed model are presented. A comparison with results predicted by the classical $M/G/c/c$ state dependent queuing model shows that
the proposed model correctly captures the interaction between upstream traffic demand and downstream traffic supply.
Future work shall derive more analytic results on the proposed model, and
include the extension of the model to complicated tree-topologies (complex series, merge,
and split networks).
Another interesting extension is the treatment of the case where the arrival rate follow a general distribution (general state dependent queuing model).
The case where traffic demand, traffic supply and fundamental diagrams are stochastic will also be treated.

\section*{Acknowledgements}
The authors wish to thank the editor-in-chief, and the anonymous reviewers whose constructive suggestions helped improve this paper.

\appendix

\section{Proof of Theorem~\ref{thm1}}
\label{appA}

\begin{enumerate}

\item \textbf{Existence.}
We have $e: [0, \lambda] \to [0, \lambda]$.\\
\begin{itemize}
 \item $e(0)=h(0)=\lambda (1-P^{(1)}_{c_1}(\lambda,0))=\lambda\left(1-\sum_{n_2=1}^{c_2} P^{(1\mid 2)}_{(c_1\mid n_2)}(\lambda) P^{(2)}_{n_2}(0)\right)$.\\
    $P^{(2)}_{n_2}(0) = 0, \forall n_2 \neq 0$ and $P^{(2)}_{0} = 1$. \\
    Therefore, $e(0) = \lambda(1 - P^{(1\mid 2)}_{(c_1\mid 0)}(\lambda)) > 0$\\
    since $P^{(1\mid 2)}_{(c_1\mid 0)}(\lambda) < 1$, because
    $P^{(1\mid 2)}_{(0\mid 0)}(\lambda) > 0$ and $\lambda > 0$.
 \item $e(\lambda) = h(\lambda) - \lambda = \lambda(1-P^{(1)}_{c_1}(\lambda))-\lambda  = -\lambda P^{(1)}_{c_1}(\lambda) < 0$\\
     since $P^{(1)}_{c_1}(\lambda) > 0$ and $\lambda > 0$.
\end{itemize}
Using the theorem of intermediate values, $e$ is continuous from $[0, \lambda]$ into $[0, \lambda]$, and we have $e(0) > 0$ and $e(\lambda) < 0$.
We conclude that $\exists \theta \in [0, \lambda]$ such that $e(\theta)=0$.

\item \textbf{Uniqueness.}
For the uniqueness of the fixed point, it is sufficient to show that $e$ is decreasing in $[0, \lambda]$. Let us calculate
$d\;e(\theta)/d \theta$.

$$\frac{d\;e(\theta)}{d \theta} = \frac{d \left( \lambda - \lambda \sum_{n_2=1}^{c_2} P^{(1\mid 2)}_{(c_1\mid n_2)}(\lambda) P^{(2)}_{n_2}(\theta) - \theta \right)}{d \theta}.$$

One can easily show that
$$\frac{d P^{(2)}_{n_2}}{d \theta} = \frac{1}{\theta} P^{(2)}_{n_2}(\theta) \left( n_2 - \bar{n}_2 (\theta)\right),$$
where
$$\bar{n}_2 (\theta) = \mathbb E_{P^{(2)}} (n_2) = \sum_{k_2=1}^{c_2} k_2 P^{(2)}_{k_2}(\theta).$$
Then
$$\frac{d\;e(\theta)}{d \theta} = - \frac{\lambda}{\theta} \sum_{n_2=1}^{c_2} P^{(1\mid 2)}_{c_1\mid n_2} (\lambda) P^{(2)}_{n_2}(n_2 - \bar{n}_2(\theta)) \; - 1.$$
Let us show that
\begin{equation}\label{eqqS}
   S \stackrel{def}{=} \sum_{n_2=1}^{c_2} P^{(1\mid 2)}_{c_1\mid n_2} (\lambda) P^{(2)}_{n_2}(\theta)(n_2 - \bar{n}_2(\theta)) \geq 0.
\end{equation}
We have $P^{(1\mid 2)}_{c_1\mid n_2} (\lambda)$ is increasing with respect to $n_2$; see~(\ref{f12}) and~(\ref{eqq12}).

Let $n^*_2 = n^*_2 (\theta)$ be defined as follows $n^*_2 = \max\{n_2 \leq c_2, n_2 - \bar{n}_2(\theta) < 0\}$.
Then
$$\begin{array}{ll}
   S   & = \sum_{n_2=1}^{n^*_2} P^{(1\mid 2)}_{c_1\mid n_2} (\lambda) P^{(2)}_{n_2}(\theta)(n_2 - \bar{n}_2(\theta))
           + \sum_{n_2=n^*_2}^{c_2} P^{(1\mid 2)}_{c_1\mid n_2} (\lambda) P^{(2)}_{n_2}(\theta)(n_2 - \bar{n}_2(\theta)) \\ \\
       & \geq \sum_{n_2=1}^{n^*_2} P^{(1\mid 2)}_{c_1\mid n^*_2} (\lambda) P^{(2)}_{n_2}(\theta)(n_2 - \bar{n}_2(\theta))
           + \sum_{n_2=n^*_2}^{c_2} P^{(1\mid 2)}_{c_1\mid n^*_2} (\lambda) P^{(2)}_{n_2}(\theta)(n_2 - \bar{n}_2(\theta)) \\ \\
       & = P^{(1\mid 2)}_{c_1\mid n^*_2} (\lambda) \sum_{n_2=1}^{c_2} P^{(2)}_{n_2}(\theta)(n_2 - \bar{n}_2(\theta)) = 0.
\end{array}$$
Hence $d\;e(\theta)/d \theta < 0$, and the fixed point is unique.

\end{enumerate}

\endproof

\section{Proof of Theorem~\ref{thm2}}
\label{appB}

Following the proof of Theorem~\ref{thm1} (~\ref{appA}), we have
$$\frac{d\; h(\theta)}{d\theta} = \frac{d\; e(\theta)}{d\theta} + 1 = -\frac{\lambda}{\theta} S,$$
where $S$ is given by~(\ref{eqqS}) in the proof of Theorem~\ref{thm1}.

The condition of Theorem~\ref{thm2} can be written $\exists 0\leq \varepsilon < \theta / \lambda$ such that $S \leq \varepsilon$.
Therefore, if this condition is satisfied, then $\exists 0\leq \eta < 1$ such that $d\;h(\theta) / d\theta > -1$.

Then since we have already shown in the proof of Theorem~\ref{thm1} (~\ref{appA}) that $s \geq 0$, and by that,
$d\;h(\theta) / d\theta \leq 0$, then, under the condition of Theorem~\ref{thm2}, we have $0 \leq d\;h(\theta) / d\theta < 1$.

Hence the fixed point iteration converges to the unique fixed point of the fixed point equation.

\endproof

\end{document}